\documentclass{amsart}

\RequirePackage{amsmath}
\RequirePackage{amsopn}
\RequirePackage{amsfonts}

\newtheorem{thm}{Theorem}[section]

\newtheorem{cor}[thm]{Corollary}
\newtheorem{lemma}[thm]{Lemma}
\newtheorem{prop}[thm]{Proposition}

\theoremstyle{definition}
\newtheorem{definition}[thm]{Definition}

\theoremstyle{remark}
\newtheorem{remark}[thm]{Remark}
\newtheorem{example}[thm]{Example}

\numberwithin{equation}{section}

\def\mathcs{C^{*}}
\def\cs{\ifmmode\mathcs\else$\mathcs$\fi}

\DeclareMathSymbol{\rtimes}{\mathbin}{AMSb}{"6F}

\def\C{\mathbf{C}}
\def\T{\mathbf{T}}

\def\K{\mathcal{K}}
\let\ker\relax

\DeclareMathOperator{\Ad}{Ad}
\DeclareMathOperator{\Ind}{Ind}

\DeclareMathOperator{\Prim}{Prim}

\DeclareMathOperator{\ker}{ker}

\DeclareMathOperator{\Aut}{Aut}

\DeclareMathOperator{\id}{id}
\def\set#1{\{\,#1\,\}}

\newbox\hidebox
\def\spechide#1{\setbox\hidebox=\hbox{$#1$}
\hbox to\wd\hidebox{$\box\hidebox^\wedge$\hss}}
\makeatletter
\def\labelenumi{\textnormal{(\@alph\c@enumi)}}
\def\theenumi{\@alph \c@enumi}
\newcount\charno
\def\alphapart#1{\charno=96
\advance\charno by#1\char\charno}
\def\partref#1{\textup{(}\textnormal{\alphapart{#1}}\textup{)}}
\makeatother
\def\<{\langle}
\def\>{\rangle}
\let\ipscriptstyle=\scriptscriptstyle
\def\lipsqueeze{{\mskip -3.0mu}}
\def\ripsqueeze{{\mskip -3.0mu}}
\def\ipcomma{\nobreak\mathrel{,}\nobreak}
\newbox\ipstrutbox
\setbox\ipstrutbox=\hbox{\vrule height8.5pt
width 0pt}
\def\ipstrut{\copy\ipstrutbox}
\def\lip#1<#2,#3>{\mathopen{\relax_{\ipstrut\ipscriptstyle{
#1}}\lipsqueeze
\langle} #2\ipcomma #3 \rangle}
\def\blip#1<#2,#3>{\mathopen{\relax_{\ipstrut
\ipscriptstyle{ #1}}\lipsqueeze\bigl\langle} #2\ipcomma #3 \bigr\rangle}
\def\rip#1<#2,#3>{\langle #2\ipcomma #3
\rangle_{\ripsqueeze\ipstrut\ipscriptstyle{#1}}}
\def\brip#1<#2,#3>{\bigl\langle #2\ipcomma #3
\bigr\rangle_{\ripsqueeze\ipstrut\ipscriptstyle{#1}}}
\def\angsqueeze{\mskip -6mu}
\def\smangsqueeze{\mskip -3.7mu}
\def\trip#1<#2,#3>{\langle\smangsqueeze\langle #2\ipcomma #3
\rangle\smangsqueeze\rangle_{\ripsqueeze\ipstrut\ipscriptstyle{#1}}}
\def\btrip#1<#2,#3>{\bigl\langle\angsqueeze\bigl\langle #2\ipcomma
#3
\bigr\rangle
\angsqueeze\bigr\rangle_{\ripsqueeze\ipstrut\ipscriptstyle{#1}}}
\def\tlip#1<#2,#3>{\mathopen{\relax_{\ipstrut\ipscriptstyle{
#1}}\lipsqueeze \langle\smangsqueeze\langle} #2\ipcomma #3
\rangle\smangsqueeze\rangle}
\def\btlip#1<#2,#3>{\mathopen{\relax_{\ipstrut\ipscriptstyle{
#1}}\lipsqueeze
\bigl\langle\angsqueeze\bigl\langle} #2\ipcomma #3
\bigr\rangle\angsqueeze\bigr\rangle}

\def\ip(#1|#2){(#1\mid #2)}
\def\bip(#1|#2){\bigl(#1 \mid #2\bigr)}
\def\Bip(#1|#2){\Bigl( #1 \bigm| #2 \Bigr)}
\def\X{{\modulefont X}}

\def\hm{homomorphism}
\def\H{\mathcal{H}}

\def\sme{\,\mathord{\mathop{\text{--}}\nolimits_{\relax}}\,}

\def\tg{\operatorname{tg}}

\def\RR{\mathbf{R}}
\def\CC{\mathbf{C}}
\def\TT{\mathbf{T}}
\def\ZZ{\mathbf{Z}}
\def\NN{\mathbf{N}}

\def\res{\operatorname{res}}

\def\uC{\underline{C}}
\def\uB{\underline{B}}
\def\uZ{\underline{Z}}
\def\uH{\underline{H}}
\def\upartial{\underline{\partial}}
\def\ab{\mathrm{ab}}
\def\loc{\mathrm{loc}}
\def\pt{\mathrm{pt}}

\def\X{\mathsf{X}}
\def\eps{\varepsilon}
\def\CR{\mathcal{CR}}

\def\zmodule_#1{\mathcal{E}_{#1}}
\def\zbarmodule_#1{\zmodule_{\underline{#1}}}
\def\dualact#1{#1^{\wedge}}
\def\idhat#1{\hat{#1}}
\def\hatu{\idhat u}
\def\hatv{\idhat v}
\def\hats{\idhat \sigma}
\def\ahat{\hat A}
\def\bhat{\hat B}
\def\charfcn_#1{\mathbb{I}_{#1}}
\def\Lprime{L'}

\usepackage[all]{xy}
\CompileMatrices
\SelectTips{cm}{}

\allowdisplaybreaks[2]

\date{15 May 2001; revised 1 March 2002}
\begin{document}

\title[Central twisted transformation groups]
{Central twisted transformation groups and
\boldmath group $C^*$-algebras of central group extensions}

\author[Echterhoff]{Siegfried Echterhoff}
\address{Westf\"alische Wilhelms-Universit\"at M\"unster \\
Mathematisches Institut \\
Einsteinstr. 62 \\
W-48149 M\"unster\\
Germany}
\email{echters@math.uni-muenster.de}

\author[Williams]{Dana P. Williams}
\address{Department of Mathematics \\
Dartmouth College \\
Hanover, NH 03755-3551 \\
USA}
\email{dana.williams@dartmouth.edu}

\begin{abstract}
We examine the structure of central twisted transformation group
\cs-algebras $C_{0}(X)\rtimes_{\id,u}G$, and apply our results to the
group \cs-algebras of central group extensions.  Our methods require
that we study Moore's cohomology group $H^{2}\bigl(G,C(X,\T)\bigr)$,
and, in particular, we prove an inflation result for pointwise trivial
cocyles which may be of use elsewhere.
\end{abstract}
\maketitle

The study of twisted transformation group \cs-algebras and Moore's
cohomology group $H^{2}\bigl(G,C(X,\T)\bigr)$ arise naturally in the
study of the group \cs-algebras of group extensions and have far
reaching applications in operator algebras.  Packer's survey articles
\cite{judy96}, and especially \cite{judy-exp}, give an excellent
historical context as well as providing additional references and the
basic definitions for what follows.

In this paper we examine the structure of central twisted transformation group
\cs-algebras $C_{0}(X)\rtimes_{\id,u}G$ where $G$ is a second
countable locally compact group acting trivially on a second countable
locally compact space $X$ and $u\in Z^{2}\bigl(G,C(X,\T)\bigr)$ is a
$2$-cocycle on $G$ taking values in the center $C(X,\T)$ of the
multiplier algebra of $C_{0}(X)$.  A primary motivation is
the study of the group \cs-algebra $\cs(L)$ of a central group
extension
\begin{equation}
    \label{eq:5}\tag{$*$}
    \xymatrix{1\ar[r]& N\ar[r] & L\ar[r]^{p}& G \ar[r] & 1.}
\end{equation}
It is well-known that $\cs(L)$ is isomorphic to the central twisted
transformation group $C_{0}(\widehat {N}) \rtimes_{\id,\eta}G$, where
$\eta \in Z^{2}(G,N)$ is a $2$-cocycle associated to the
extension~\eqref{eq:5} viewed as taking values in $C(\widehat N,\T)$.

Our results apply, in particular, to groups $G$ which are smooth in
the sense of Moore \cite{moore2} in that there is a central group
extension
\begin{equation*}
    \xymatrix{1\ar[r]& Z\ar[r] & H\ar[r]^{q}& G \ar[r] & 1}
\end{equation*}
(called a representation group for $G$) such that the transgression
map $\tg:\widehat Z\to H^{2}(G,\T)$ is an isomorphism of topological
groups. Thus our results apply in particular to groups $G$ that are
either discrete, compact, compactly
generated abelian or connected, simply connected Lie groups.
Our results substantially generalize results
in \cite{hasm1, hasm2, 90a}, where central twisted group
algebras have been studied for abelian groups
and pointwise trivial cocycles $u\in Z^2(G, C(X,\TT))$,
and results in \cite{hasm3}, where Smith considers
central twisted group algebras for groups
with $H^2(G,\TT)$ discrete (which, for example, is the case
whenever $G$ is compact).

The structure of $C_{0}(X)\rtimes_{\id,u}G$ is easiest to describe
when $u$ is pointwise trivial; that is, $u(x)$ is trivial in
$H^{2}(G,\T)$ for each $x\in X$, where $u(x)(s,t):=u(s,t)(x)$.  As was
shown in \cite{hasm1, 90a}
for abelian $G$, each pointwise trivial cocycle is
naturally associated to a locally compact $\widehat{G}_{\ab}$-bundle
$p:\zmodule_{u} \to X$, where $G_{\ab}:=G/\overline{[G,G]}$ is the
\emph{abelianization} of $G$ (Definition~\ref{def-Zu}).  The bundles
$\zmodule_{u}$ and $\zmodule_{v}$ are isomorphic if and only if $u$
and $v$ are equivalent in $H^{2}\bigl(G,C(X,\T)\bigr)$.  Furthermore,
the map $u\mapsto \zmodule_{u}$ is multiplicative in that
$\zmodule_{uv}\cong \zmodule_{u}*\zmodule_{v}$, where the
$\zmodule_{u}*\zmodule_{v}$ is the usual product of
$\widehat{G}_{\ab}$-bundles. In the case of abelian $G$, the 
bundles $\zmodule_{u}$ appeared in \cite{hasm2} under the name of
{\em characteristic bundles} the
isomorphism classes of which form an (often proper) subset
of the set of all isomorphism classes of free and proper
$\widehat{G}$-bundles over $X$.
Moreover, in this case, the algebra
$C_{0}(X)\rtimes_{\id,u}G$ is isomorphic to
$C_{0}(\zmodule_{u})$. We shall give a short proof
of this well-known isomorphism below.

When $G$ is not abelian,
$C_{0}(X)\rtimes_{\id,u}G$ is still a $C_{0}(X)$-algebra, and
therefore
can be
thought of as the (semi-continuous) sections of a \cs-bundle over $X$.
If $A$ is a $C_{0}(X)$-algebra admitting a $C_{0}(X)$-linear
$\widehat{G}_{\ab}$-action, then we can form the \cs-analogue of the
bundle product above, and obtain a \cs-algebra $\zmodule_{u}*A$
(\cite[Definition~3.3]{ew3}).  One benefit of this construction is
that
if $u$ is
pointwise trivial and $G$ is smooth then we can show that
\begin{equation}
    \label{eq:3}\tag{$\dag$}
    C_{0}(X)\rtimes_{\id,u}G\cong
    \zmodule_{u}*C_{0}\bigl(X,\cs(G)\bigr)
\end{equation}
(Theorem~\ref{thm-pt}).  A crucial ingredient in the proof of
\eqref{eq:3}, which may be of independent interest, is showing that
$u$ is equivalent to a cocycle inflated from a pointwise trivial
cocycle on $G_{\ab}$ (Proposition~\ref{prop-lifting}).  More
generally, we show that
\begin{equation}
    \label{eq:7}\tag{$\ddag$}
    H^{2}\bigl(G,C(X,\T)\bigr) \cong C\bigl(X,H^{2}(G,\T)\bigr) \oplus
    H_{\pt}^{2} \bigl(G_{\ab},C(X,\T)\bigr)
\end{equation}
(Proposition~\ref{prop-smooth}).
If $u\in Z^{2}\bigl(G,C(X,\T)\bigr)$, then the map $x\mapsto [u(x)]$
gives us an element $\varphi\in C\bigl(X,H^{2}(G,\T)\bigr)$ and we can
use \eqref{eq:7} to write $u=v\cdot {\bar u}_{\varphi}$ with $v$
pointwise trivial.  Since $H$ is a representation group for $G$,
$\cs(H)$ is a $C_{0}(\widehat Z)$-algebra.  Since we can view $\varphi$ as
a continuous map of $X$ into $\widehat Z\cong H^{2}(G,\T)$,
we can form the pull-back $\varphi^{*}\bigl(\cs(H)\bigr)$.
Our main result (Theorem~\ref{thm-general}) implies that
\begin{equation*}
    C_{0}(X)\rtimes_{\id,u}G\cong
    \zmodule_{v}*\varphi^{*}\bigl(\cs(H)\bigr) .
\end{equation*}
This result is new even in the special situation
  where $H^2(G,\TT)$ is discrete as is the case in \cite{hasm3}
(see Theorem \ref{thm-pt} below).

In the case of a central group extension $\cs(L)\cong C_{0}(\widehat
N)\rtimes _{\id,\eta}G$, our results take a more elegant form as our
auxiliary constructions can be formulated in group theoretic terms.
For example, if $L$ is pointwise trivial --- that is, the associated
cocycle $\eta$ is pointwise trivial --- then
$\res:\widehat{L}_{\ab}\to \widehat N$ is a $\widehat{G}_{\ab}$-bundle
which is isomorphic to $\zmodule_{\eta}$, and we have
\begin{equation*}
    C^{*}(L)\cong \widehat{L}_{\ab}*C_{0}\bigl(\widehat N, C^{*}(G)\bigr)
\end{equation*}
(Corollary~\ref{cor-pt-group}).
In general, we prove (Theorem~\ref{thm-groupex}) that
\begin{equation*}
    C^{*}(L)\cong
    \widehat{\Lprime}_{\ab}*\varphi^{*}\bigl(C^{*}(H)\bigr),
\end{equation*}
where $\Lprime$ is the quotient of $\set{(l,h)\in L\times H:
    p(l)=q(h)}$ by the subgroup $\Delta(Z):=\set{(\hat\varphi(z),z):z\in
    Z}$ where $\hat\varphi:Z\to N$ is the dual to $\varphi:\widehat N\to
    H^{2}(G, \T)\cong\widehat Z$.

The paper is organized as follows.  In Section~\ref{sec:pointw}, we
make a careful study of pointwise trivial cocycles and prove our
results on
inflation of pointwise trivial cocycles, and on the decomposition of
$H^{2}\bigl(G,C(X,\T)\bigr)$. In
Section~\ref{sec:centr}, we prove our structure theorems for general
central twisted transformation groups.  In Section~\ref{sec:group}, we
show how to apply our results to the group \cs-algebras of central
group extensions.

We should mention that this paper is very much related to the papers
\cite{ew2, ew3}, where we consider more general systems.
However, in the special situations considered here,
the results are much easier to describe, and allow more general
statements.

\section{Pointwise Trivial Cocycles}\label{sec:pointw}

If $X$ is a second countable locally compact space, then the set of
continuous functions $C(X,\T)$ from $X$ into the circle group $\T$ is
a Polish group when equipped with pointwise multiplication and the
compact-open topology.  In this section, we want to look carefully at
the Moore cohomology group $H^2\bigl(G, C(X,\TT)\bigr)$ for a second
countable locally compact group $G$ acting trivially on $C(X,\T)$.

The definition and basic properties of Moore's cohomology groups
$H^{n}(G,A)$ for an arbitrary polish $G$-module $A$ are laid out in
Moore's original paper \cite{moore3}.  (A summary with additional
references can be found in \S7.4 of \cite{rw-book}.)  An important
facet of the theory is that these groups can be computed by using
two different complexes: one can either take the complex
\begin{equation*}
    \xymatrix@1{\cdots\ar[r]^-{\partial}& C^n(G,A) \ar[r]^-{\partial}&
   C^{n+1}(G,A) \ar[r]^-{\partial}&\cdots,}
\end{equation*}
where $C^n(G,A)$ denotes the group of $A$-valued
Borel functions on $G^n$ and
$\partial$ denotes the usual group coboundary, or one can work with
the complex
\begin{equation*}
    \xymatrix@1{\cdots\ar[r]^-{\protect{\upartial}}& \protect{\uC}^n(G,A)
\ar[r]^-{\protect{\upartial}}&
   \protect{\uC}^{n+1}(G,A) \ar[r]^-{\protect{\upartial}}&\cdots,}
\end{equation*}
where $\uC^{n}(G,A)$ is the quotient of $C^n(G,A)$ obtained by
identifying Borel functions on $G^n$ which coincide almost everywhere
and $\upartial$ is the induced map.  Moore shows in
\cite[Theorem~5]{moore3} that the canonical maps $C^n(G,A)\to
\uC^n(G,A)$ induce isomorphisms of $H^n(G,A)$ with $\uH^n(G,A)$ for
all $n\ge0$.  One advantage of working with $\uC^n(G,A)$ is that
$\uC^n(G,A)$ is a Polish group when equipped with the topology of
convergence in measure (after replacing Haar measure with an
equivalent finite measure).  Thus we have a topology on $\uH^n(G,A)$
(and therefore on $H^n(G,A)$), although this topology can be
non-Hausdorff in general.  On the other hand, elements in $\uC^n(G,A)$
are not defined everywhere, and this can often be a nuisance; thus it
is useful to work with both definitions simultaneously, and we shall
do so below.

\begin{definition}\label{def-pointwise}
    Let $u\in Z^2(G, C(X,\TT))$ and define $u(x)\in Z^2(G,\TT)$ by
    evaluation at $x$: $u(x)(s,t):=u(s,t)(x)$.  We say that $u$ is
    \emph{pointwise trivial} if $u(x)$ is trivial (i.e., $u(x)\in
    B^2(G,\TT)$) for all $x\in X$.  We say that $u$ is \emph{locally
      trivial} if each $x\in X$ has an open neighborhood $V$ such that
    the restriction $u_V\in Z^2\bigl(G, C(V,\TT)\bigr)$ of $u$ to $V$ is
    trivial (i.e., $u_V\in B^2\bigl(G, C(V,\TT)\bigr)$).
\end{definition}

We denote by $Z^2_{\pt}(G, C(X,\TT))$ and $Z^2_{\loc}(G, C(X,\TT))$,
respectively, the subsets of pointwise trivial cocycles and locally
trivial cocycles in $Z^2(G, C(X,\TT))$.  Similarly, we let
$H^2_{\pt}(G, C(X,\TT))$ and $H^2_{\loc}(G, C(X,\TT))$ be the images
of $Z^2_{\pt}(G, C(X,\TT))$ and $Z^2_{\loc}(G, C(X,\TT))$ in $H^2(G,
C(X,\TT))$, respectively.

Pointwise trivial cocycles have been studied extensively in the
literature. A particularly important example --- having applications
in the study of $C^*$-dynamical systems --- is Rosenberg's
\cite[Theorem 2.1]{ros2}, which shows that $Z^2_{\pt}\bigl(G,
C(X,\TT)\bigr)=Z^2_{\loc}\bigl(G, C(X,\TT)\bigr)$, whenever
$H^2(G,\TT)$ is Hausdorff and the abelianization
$G_{\ab}$ is compactly generated.  However,
Rosenberg's Theorem can fail without the assumptions on $G_{\ab}$ and
$H^2(G,\TT)$ (see \cite{ew2} and Example~\ref{ex-open} below).

We shall study pointwise trivial cocycles $u\in Z^2_{\pt}\bigl(G,
C(X,\TT)\bigr)$ via a canonical $\widehat{G}_{\ab}$-space
$\zmodule_u$.  Our construction of $\zmodule_u$ is identical to that
in \cite{90a} where the group $G$ was assumed to be abelian. It
follows from \cite[Theorem~3]{moore3} that $H^{1}(G,\T)$, and
therefore $\uH^{1}(G,\T)$, can be identified with
$\widehat{G}_{\ab}$.  Since there are natural maps of
$\widehat{G}_{\ab}$ into $C^{1}(G,\T)$ and $\uC^{1}(G,\T)$, we always
have \emph{algebraic} short exact sequences of groups
\begin{gather*}
    \xymatrixnocompile@1{1 \ar[r] &\widehat{G}_{{\ab}} \ar[r] &
      C^{1}(G,{\T})
      \ar[r]^{\partial} & B^{2}(G,{\T}) \ar[r] & 1}\\
    \intertext{and} \xymatrixnocompile@1{1 \ar[r] &\widehat{G}_{{\ab}}
      \ar[r] & {\uC}^{1}(G,{\T}) \ar[r]^{{\upartial}} &
      {\uB}^{2}(G,{\T}) \ar[r] &1, }
\end{gather*}
which are related to each other via the inclusions $C^n\to \uC^n$.
Although the homomorphisms in the second sequence are always
continuous, the second sequence may fail to be a short exact sequence
of topological groups as the associated quotient map
$\uC^1(G,\TT)/\widehat{G}_{\ab}\to \uB^2(G,\TT)$ may fail to be a
homeomorphism.  In fact, this quotient map is a topological
isomorphism if and only if $\uB^2(G,\TT)$ is a Polish subgroup of
$\uC^2(G,\TT)$.  This happens exactly when $\uB^2(G,\TT)$ is closed
which is equivalent to $\uH^2(G,\TT)$ being Hausdorff (cf., e.g.,
\cite[\S5]{moore3}).

\begin{definition}\label{def-Zu}
    Let $u\in Z^2_{\pt}\bigl(G, C(X,\TT)\bigr)$. Then we define
\begin{equation*}
\zmodule_u=\set{(f, x)\in C^1(G,\TT)\times X:
       \partial(f)=u(x)}.
\end{equation*}
Similarly, if $\underline{u}$ denotes the image of $u$ in
$\uZ^2\bigl(G, C(X,\TT)\bigr)$, we define
\begin{equation*}
\zbarmodule_u =\set{(\underline{f},x)\in \uC^1(G,\TT)\times X:
       \upartial(\underline{f})=\underline{u}(x)}.
\end{equation*}
\end{definition}

It follows from the short exact sequences mentioned above that both,
$\zmodule_u$ and $\zbarmodule_u $ are free $\widehat{G}_{\ab}$-spaces,
and that $\zbarmodule_u $ becomes a topological
$\widehat{G}_{\ab}$-space when equipped with the relative topology
from $\uC^1(G,\TT)\times X$. Moreover, the canonical projections
\begin{equation*}
p:\zmodule_u\to X\quad\text{and}\quad \underline{p}:\zbarmodule_u
\to X\end{equation*}
induce bijections between the quotient spaces
$\zmodule_{u}/\widehat{G}_{\ab}$ and $\zbarmodule_u
/\widehat{G}_{\ab}$ and $X$, respectively. However,  there are cases where
$\underline{p}:\zbarmodule_u \to X$ fails
to be open (e.g., Example~\ref{ex-open} below).

\begin{prop}\label{prop-bijection} Let $u\in Z^2_{\pt}\bigl(G, C(X,\TT)\bigr)$.
    Then the following are true:
       \begin{enumerate}
       \item The map $ (f, x)\mapsto (\underline{f},x)$ is a bijection
         $\varphi:\zmodule_u\to \zbarmodule_u$.
       \item If we topologize $\zmodule_u$ via the identification with
         $\zbarmodule_u $ of \partref1, then the topology on
         $\zmodule_u\subset C^1(G,\TT)\times X$ is given by pointwise
         convergence in the first variable, and the given topology on
         $X$.
       \item If $[u]=[v]\in H^2_{\pt}\bigl(G, C(X,\TT)\bigr)$, then
         $\zmodule_u\cong \zmodule_v$ as (topological)
         $\widehat{G}_{\ab}$-spaces.
       \item $\zmodule_u$ is isomorphic to the trivial
         $\widehat{G}_{\ab}$-bundle $\widehat{G}_{\ab}\times X$ (as a
         $\widehat{G}_{\ab}$-space) if and only if $u$ is trivial.
      \end{enumerate}
\end{prop}
\begin{proof}
    For the proof of \partref1, it is enough to show that the given map
    induces bijections between the $\widehat{G}_{\ab}$-orbits in
    $\zmodule_u$ and $\zbarmodule_u $. But, by the above discussion, the
    orbits are just the sets $p^{-1}(\set{x})$ and
    $\underline{p}^{-1}(\set{x})$, respectively, and it is clear that
    $p^{-1}(\set{x})$ is mapped into $\underline{p}^{-1}(\set{x})$. The
    result then follows from the freeness of the
    $\widehat{G}_{\ab}$-actions.

    For \partref2, we have to show that a sequence $\set{(f_n, x_n)}$
    converges to $(f, x)$ in $\zmodule_u$ if and only if $f_n\to f$
    pointwise on $G$ and $x_n\to x$.  Since pointwise convergence of
    $\set{f_n}$ in $C^1(G,\TT)$ implies convergence of
    $\set{\underline{f_n}}$ in $\uC^{1}(G,\T)$, it is enough to show
    that convergence of $\set{(\underline{f_n}, x_n)}$ in $\zbarmodule_u
    $ implies pointwise convergence of $\set{f_n}$.  But this follows
    from part~\partref1 and \cite[Lemma 3.6]{90a}.

    For \partref3, let $g\in C^1\bigl(G, C(X,\TT)\bigr)$ be a Borel
    cochain such that $v=\partial(g)u$. Then $\zmodule_u\to \zmodule_v:
    (f,x)\mapsto (g(x)f, x)$ is an isomorphism which is bicontinuous by
    part~\partref2.

    The last assertion follows from the proof of \cite[Proposition
    3.4]{90a}, which did not make use of the assumption that $G$ was
    abelian.
\end{proof}

\begin{remark}\label{rem-local}
    It follows from part~\partref4 of the proposition that $\zmodule_u$
    is a locally trivial $\widehat{G}_{\ab}$-bundle if and only if $u$
    is locally trivial. In particular, when $u$ is locally trivial,
    $\zmodule_u$ is a free and proper locally compact
    $\widehat{G}_{\ab}$-bundle.
\end{remark}

If $G$ is abelian, then \cite[Proposition 3.1]{90a} implies that
$\zmodule_u$ is locally compact, and $p:\zmodule_u\to X$ is a free and
proper $\widehat{G}$-bundle. A careful look at the proof reveals that
the hypothesis that $G$ be abelian was only used to guarantee that the
coboundary map $\upartial:\uC^1(G,\TT)\to \uB^2(G,\TT)$ is
open.\footnote{If
$G$ is abelian, then $H^{2}(G,\T)$ is Hausdorff
\cite[Theorem~7]{moore4}, and this is equivalent to the openness of
$\upartial$.}  Thus we get

\begin{prop}\label{prop-proper}
    Assume that $H^2(G,\TT)$ is Hausdorff and $u\in Z^2_{\pt}(G,
    C(X,\TT))$. Then $\zmodule_u$ is locally compact and
    $p:\zmodule_u\to X$ is a free and proper $\widehat{G}_{\ab}$-bundle.
\end{prop}

We shall also need the following lemma which can be proved along the
lines of the final part of the proof of \cite[Proposition 3.4]{90a}.

\begin{lemma}\label{lem-borel}
    Let $Y$ be a second countable topological space and let $g:Y\to
    C(X,\TT)$ be a map such that for each $x\in X$ the map $g(x):Y\to
    \TT: g(x)(y)=g(y,x)$ is Borel.  Then $g:Y\to C(X,\T)$ is Borel.
\end{lemma}

We are now ready for our lifting result for pointwise trivial
cocycles.

\begin{prop}\label{prop-lifting}
    Assume that $u\in Z^2_{\pt}\bigl(G, C(X,\TT)\bigr)$. Then the
    following are equivalent:
       \begin{enumerate}
       \item There exists a cocycle $\tilde{u}\in Z^2_{\pt}(G_{\ab},
         C(X,\TT))$ such that $u$ is cohomologous to the inflation
         $\inf\tilde{u}$ of $\tilde{u}$ to $G$.
       \item The projection $p:\zmodule_u\to X$ is open.
       \item $\zmodule_u$ is locally compact and $p:\zmodule_u\to X$ is
         a free and proper $\widehat{G}_{\ab}$-bundle.
        \end{enumerate}
        Moreover, if \partref1 holds, then $\zmodule_u$ is isomorphic to
        $\zmodule_{\tilde{u}}$ as a $\widehat{G}_{\ab}$-space.
\end{prop}
\begin{proof} We will first show the last statement:
    indeed, if $u\sim \inf\tilde{u}$, we can use part~\partref3 of
    Proposition~\ref{prop-bijection} to assume without loss of
    generality that $u=\inf\tilde{u}$.  It is then easy to check that
    the map from $\zmodule_{\tilde{u}}\to \zmodule_u$ given by
    \begin{equation*}
    (f, x)\mapsto (\inf f, x)
    \end{equation*}
    is a $\widehat{G}_{\ab}$-equivariant bijection, which is
    bicontinuous by part~\partref2 of Proposition~\ref{prop-bijection}.

    We now show \partref1 $\Longrightarrow$~\partref3
    $\Longrightarrow$~\partref2~$\Longrightarrow$~\partref1.  In fact
    \partref1~$\Longrightarrow$~\partref3 follows from the isomorphism
    $\zmodule_u\cong \zmodule_{\tilde{u}}$ together with
    \cite[Proposition 3.1]{90a} (see the discussion preceding
    Proposition \ref{prop-proper} above), and
    \partref3~$\Longrightarrow$~\partref2 follows by definition.

    We now check \partref2~$\Longrightarrow$~\partref1.  Choose a Borel
    section $c:G_{\ab}\to G$ with $c(\dot{e})=e$ and define a map $\mu:
    G_{\ab}\times G_{\ab}\times \zmodule_u\to\TT$ by
    \begin{equation*}
   \mu\big(\dot{s},\dot{t}, (f,x)\big)= \partial_{G_{\ab}}(f\circ
    c)(\dot{s},\dot{t}).\end{equation*}
    This map is continuous in $(f,x)$ by
    Proposition~\ref{prop-bijection}\partref2, and since
    $\partial_{G_{\ab}}\big(\gamma\cdot (f\circ c)\big)=
    \partial_{G_{\ab}}(f\circ c)$ for all $\gamma\in \widehat{G}_{\ab}$,
    it follows that $\mu$ induces a map $\tilde{u}:G_{\ab}\times
    G_{\ab}\times X\to \TT$ which is continuous on $X$ (since
    part~\partref2 implies that $p:\zmodule_u\to X$ induces a
    homeomorphism
    $\zmodule_u/\widehat{G}_{\ab}\cong X$).  Thus we may view
    $\tilde{u}$ as a map from $G_{\ab}\times G_{\ab}$ to $C(X,\T)$. Since
    for each $x\in X$, $\tilde u(x)=\partial_{G_{\ab}}(f\circ c)$ for
    any
    $(f,x)\in
    \zmodule_u$, it follows that all evaluations $\tilde{u}(x)$ are
    Borel.  Therefore, $\tilde{u}$ is Borel by Lemma~\ref{lem-borel}.

    It remains to check that $u\sim \inf\tilde{u}$.  For this we define
    a map $\nu: G\times \zmodule_u\to \TT$ by
    \begin{equation*}
\nu\bigl(s, (f,x)\bigr)=\overline{f(s)}\cdot \big(f\circ
c\big)(\dot{s}).\end{equation*}
    Then $\nu$
    is continuous on $\zmodule_u$, and since $\nu\bigl(s, (\gamma\cdot
    f,x)\bigr)=\nu\bigl(s, (f,x)\bigr)$ for all $\gamma\in 
\widehat{G}_{\ab}$, it
    follows that $\nu$ factors through a map $g: G\times X\to \TT$ which
    is continuous on $X\cong Z_{u}/\widehat{G}_{\ab}$.
    Lemma~\ref{lem-borel} implies that $g$, viewed as a map from $G$ to
    $C(X,\TT)$, is Borel. Moreover, if $(f,x)\in \zmodule_u$, then
    \begin{equation*}
\partial g(s,t)(x)=\partial \overline{f}\cdot \partial (f\circ c)
    =\overline{u(x)}\cdot \inf\tilde{u}(x).\end{equation*}
    This completes the proof.
\end{proof}

\begin{example}[{cf., \cite[Example 7.3]{ew2}}]\label{ex-open}
    There exist groups $G$ with cocycles $u\in Z^2_{\pt}\bigl(G,
    C(X,\TT)\bigr)$ such that $p:\zmodule_u\to X$ fails to be open; it
    follows that such $u$ are not lifted from pointwise trivial cocycles
    on $G_{\ab}$. Let $\theta$ be an irrational number and let
    $G=\RR^2\times \TT^2$ with multiplication given by the formula
    \begin{equation*}
(s_1, t_1, z_1, w_1)(s_2, t_2, z_2, w_2) =(s_1+s_2, t_1+t_2,
    e^{is_1t_2}z_1z_2, e^{i\theta s_1 t_2}w_1w_2).\end{equation*}
    This is the
    example of a group with non-Hausdorff $H^2(G,\TT)$ presented by
    Moore in \cite[p.~85]{moore2} (see also \cite[Example 7.2]{ew2}).

    Let $X=\set{\frac{1}{n}:n\in \NN}\cup \set{0}$ and choose a
    continuous map $\lambda: X\to \ZZ+\theta\ZZ\subseteq \RR$ such that
    $\lambda_0=0$ and such that $\lambda_{\frac{1}{n}}\neq 0$ for all
    $n\in \NN$. Define a cocycle $v\in Z^2\bigl(\RR^2, C(X,\TT)\bigr)$
    by
    \begin{equation*}
v\bigl((s_1,t_1), (s_2, t_2)\bigr)(x)= e^{i\lambda(x) s_1t_2},\quad
    s,t\in G, x\in X.\end{equation*}
    Since the $v(x)$ are non-trivial cocycles
    on $\RR^2= G_{\ab}$ if $x\neq 0$, it follows that $v$ is not
    pointwise trivial.  However, the inflation $\inf v\in Z^2\bigl(G,
    C(X,\TT)\bigr)$ \emph{is} pointwise trivial.  Indeed, a short
    computation shows that if we define
    $\lambda(x)=l(x)+\theta m(x)$ where $l,m:X\to \ZZ$, and if we define
    $f_x\in C^1(G,\TT)$ by $f_x(s,t,z,w)=z^{l(x)}w^{m(x)}$, then $\inf
    v(x)=\partial f_x$.

    We now show that the projection $p:\zmodule_{\inf v}\to X$ is not
    open.  If it were, then Proposition~\ref{prop-lifting} would imply
    that $\zmodule_{\inf v}$ would be a locally compact free and proper
    $\widehat{G}_{\ab}\cong \RR^2$-space.  Since every free and proper
    $\RR^2$-bundle is locally trivial by Palais's Slice Theorem
    \cite[Theorem~4.1]{palais}, and therefore trivial since $\check
    H^1(X,\RR^2)=\set{0}$, this would imply that $\zmodule_{\inf v}$
    would be a trivial bundle.  In that case, there exists a continuous
    section $\varphi: X\to \zmodule_{\inf v}$. Then we could find
    elements $\gamma_x\in \widehat{\RR}^{2}$ such that
    \begin{equation*}
\varphi(x)=(\gamma_x\cdot f_x, x),\quad\text{with}\quad
\gamma_0=1_G,
\end{equation*}
and where $f_x$ is defined as above.  This would imply that
$\gamma_{\frac{1}{n}}\cdot f_{\frac{1}{n}}$ converges pointwise to
$1_G$, and hence, that $f_{\frac{1}{n}}|_{\TT^2}$ converges pointwise
to $1_{\TT^2}$.  But this is impossible since for all $x\neq 0$,
$f_x|_{\TT^2}$ is a non-trivial character of $\TT^2$,
$\widehat{\TT}^{2}\cong \ZZ^2$ is discrete, and pointwise convergence
of characters implies convergence \cite[Theorem~8]{moore3}.
\end{example}

\begin{cor}\label{cor-Hausdorff}
    Assume that $H^2(G,\TT)$ is Hausdorff. Then the inflation map
    \begin{equation*}
\inf:H^2_{\pt}\bigl(G_{\ab}, C(X,\TT)\bigr)\to H^2_{\pt}\bigl(G,
    C(X,\TT)\bigr)
\end{equation*}
sending $[v]\mapsto [\inf v]$ is an isomorphism of abelian groups.

Similarly, if $G$ is any second countable locally compact group, then
the inflation map $\inf: H^2_{\loc}\bigl(G_{\ab}, C(X,\TT)\bigr)\to
H^2_{\loc}\bigl(G, C(X,\TT)\bigr)$ is an isomorphism.
\end{cor}

\begin{proof} It follows from Propositions \ref{prop-proper}~and
    \ref{prop-lifting} that
\begin{equation*}
\inf:H^2_{\pt}(G_{\ab}, C(X,\TT))\to
    H^2_{\pt}\bigl(G, C(X,\TT)\bigr)
\end{equation*}
  is surjective.  Injectivity
    follows since $[\inf v]=[1]$ in $H^2_{\pt}\bigl(G, C(X,\TT)\bigr)$
    if and only if $\zmodule_{\inf v}$ is the trivial bundle
    (Proposition~\ref{prop-bijection}\partref4).  But $\zmodule_{\inf
      v}\cong \zmodule_v$ by Proposition~\ref{prop-lifting}, and this
    implies $[v]=1$.  The second statement is proved similarly using
    Remark~\ref{rem-local} and Proposition~\ref{prop-lifting}.
\end{proof}

The above results can be used to give a description of
$H^2\bigl(G,C(X,\TT)\bigr)$ along the lines of \cite[\S5]{ew2}.  We
restrict our attention to groups which are \emph{smooth} in the sense
of Moore (see \cite{moore2} --- an extensive discussion of smooth
groups is also given in \cite[\S4]{ew2}). Recall that if $1\to Z\to
H\to G\to 1$ is a central group extension, then the
\emph{transgression map}
\begin{equation*}
\tg:\widehat{Z}\to H^2(G,\TT)
\end{equation*}
is defined by composing the characters of $Z$ with a cocycle $\eta\in
Z^2(G, Z)$ corresponding to the extension (recall that such an $\eta$
is given by $\eta=\partial c$ for any Borel section $c:G\to H$).  A
group $G$ is called smooth if there exists a central extension as
above such that $\tg:\widehat{Z}\to H^2(G,\TT)$ is bijective, which
automatically implies that it is an isomorphism of topological groups.
The extension $H$ is then called a \emph{representation group} for
$G$.  In particular, if $G$ is smooth, then $H^2(G,\TT)$ is Hausdorff.
The list of smooth groups is quite large; it contains all
discrete groups, all compact groups, all compactly generated abelian
groups, and all simply connected and connected Lie groups (see
\cite{moore3, ew2}).

Suppose now that $G$ is smooth and that $1\to Z\to H\to G\to 1$ is a
representation group of $G$. Let $\eta\in Z^2(G, Z)$ be a
corresponding cocycle.  Then any continuous map $\varphi:X\to
H^2(G,\TT)\cong \widehat{Z}$ determines a cocycle $u_{\varphi}\in
Z^2\bigl(G, C(X,\TT)\bigr)$ by defining
\begin{equation*}
u_{\varphi}(s,t)(x)=\varphi(x)\circ \eta(s,t).
\end{equation*}
It is easy to check (and it follows from the proof of
\cite[Theorem~5.4]{ew2}) that $\varphi\mapsto [u_{\varphi}]$
determines a well defined group homomorphism
\begin{equation*}
\Psi_H: C\bigl(X, H^2(G,\TT) \bigr)\to
H^2\bigl(G, C(X,\TT)\bigr),
\end{equation*}
which depends only on the particular choice of the representation group,
but not on the
particular choice of $\eta\in Z^2(G,Z)$ corresponding to this
extension.

\begin{prop}\label{prop-smooth}
    Suppose that $G$ is smooth with representation group $H$. Then the
    map
    \begin{equation*}
\inf\oplus \Psi_H:H^2_{\pt} \bigl(G_{\ab}, C(X,\TT)\bigr) \oplus C \bigl(X,
    H^2(G,\TT)\bigr)\to H^2\bigl(G, C(X,\TT)\bigr),
\end{equation*}
sending $(v, \varphi)\mapsto \inf v\cdot u_{\varphi}$, is an
isomorphism of groups.
\end{prop}
\begin{proof}
    Let $\Phi:H^2\bigl(G, C(X,\TT)\bigr)\to C\bigl(X,H^2(G,\TT)\bigr)$
    be the evaluation map given by $\Phi([u])(x):=[u(x)]$. Then $\Phi$ is a
    group homomorphism and $\ker\Phi=H^2_{\pt}\bigl(G, C(X,\TT)\bigr)$.
    Since $H^2(G,\TT)$ is Hausdorff, we can apply
    Corollary~\ref{cor-Hausdorff} to see that
    $\inf:H^2_{\pt}\bigl(G_{\ab}, C(X,\TT)\bigr)\to H^2_{\pt}\bigl(G,
    C(X,\TT)\bigr)=\ker\Phi$ is an isomorphism.  The result then follows
    from the fact that $\Psi_H:C\bigl(X, H^2(G,\TT)\bigr)\to H^2\bigl(G,
    C(X,\TT)\bigr)$ is a splitting homomorphism for $\Phi$ (see the
    proof of \cite[Theorem 5.4]{ew2} for more details).
\end{proof}

\begin{remark}\label{rem-disc}
  (1) If $H^2(G,\TT)$ is discrete, then every continuous map
  $\varphi:X\to H^2(G,\TT)$ is locally constant. Thus, together with
  the classification of characteristic bundles given in \cite{hasm2},
  the above result subsumes \cite[Corollary 9]{hasm3} --- without
  having to assume that $G$ is abelian as in \cite{hasm3}!

(2) Since every countable discrete group $G$ has a representation group
by \cite[Theorem 3.1]{moore2} (see also \cite[Corollary 1.3]{para3}),
the above decomposition applies
to such groups. If, in addition, $G_{\ab}$ is free abelian,
then it follows from the classification of characteristic bundles
in \cite[Lemma 3]{hasm2} (but see also \cite{judy96}), that
$H^2_{\pt}\bigl(G_{\ab}, C(X,\TT)\bigr)=\{0\}$.
Thus Proposition \ref{prop-smooth} implies that
\begin{equation}\label{eq:6}
H^2\bigl(G, C(X,\TT)\bigr)\cong C\bigl( X, H^2(G,\TT)\bigr)
\end{equation}
in these cases.  Notice that if $G$ is a non-abelian free group, then
$H^2(G,\TT)=\{0\}$.  Then \eqref{eq:6} implies the well know result that
$H^2\bigl(G, C(X,\TT)\bigr)=\{0\}$.
\end{remark}

\section{Central Twisted Transformation Group
   Algebras}\label{sec:centr}

In this section, we want to give a description of the central twisted
transformation group algebra $C_0(X)\rtimes_{\id,u}G$ corresponding to
a cocycle $u\in Z^2\bigl(G, C(X,\T)\bigr)$.  The basic properties of
twisted crossed products are given in \cite{para1}.  We recall some of
the fundamentals here. A (Busby-Smith) \emph{twisted action} $(\alpha,u)$ of
$G$ on a $C^*$-algebra $A$ consists of a strongly measurable map
$\alpha:G\to\Aut A $ together with a strictly measurable map
$u:G\times G\to UM(A)$ such that
\begin{equation*}
\alpha_s\alpha_t=\Ad u(s,t)\circ \alpha_{st}\quad
\text{and}\quad \alpha_r\bigl(u(s,t)\bigr)u(r,st)=u(r,s)u(rs,t),
\end{equation*}
for all $s,t,r\in G$.  We also require that $\alpha_e=\id$ and
$u(e,s)=u(s,e)=1$ for all $s\in G$.  The \emph{twisted crossed
   product} $A\rtimes_{\alpha,u}G$ is a completion of $L^1(G,A)$ with
convolution defined by
\begin{equation*}
f*g(s)=\int_G f(t)\alpha_t\bigl(g(t^{-1}s)\bigr)u(t,t^{-1}s)\,ds.
\end{equation*}
The \emph{covariant representations} of the twisted system
$(A,G,\alpha,u)$ are the pairs $(\pi,U)$ in which $\pi:A\to B(\H)$ is
a nondegenerate $*$-representation of $A$ and $U:G\to U(\H)$ is a
measurable map such that
    \begin{equation*}
U_e=1,\quad \pi\bigl(\alpha_s(a)\bigr)=U_s\pi(a)U_s^*,\quad\text{and}
       \quad
       U_sU_t=\pi\bigl(u(s,t)\bigr)U_{st}.
\end{equation*}
There is a natural one-to-one correspondence between covariant
representations $(\pi,U)$ of $(A,G,\alpha,u)$ and nondegenerate
$*$-representations of $A\rtimes_{\alpha,u}G$ associating $(\pi, U)$
to its \emph{integrated form}
\begin{equation*}
\pi\rtimes U(f)=\int_G\pi(f(s)) u_s\,ds, \quad f\in L^1(G,A).
\end{equation*}
The \emph{dual action} $\dualact{(\alpha,u)}$ of $\widehat{G}_{\ab}$
on $A\rtimes_{\alpha,u}G$ is defined on the $L^1$-functions via
\begin{equation*}
\dualact{(\alpha,u)}_{\gamma}(f)(s)=\overline{\gamma}(s)f(s).
\end{equation*}
(Just as with the Fourier transform, there is a choice to be made when
defining the dual action, and it is a matter of convenience whether
one multiplies with $\gamma(s)$ or $\overline{\gamma(s)}$ in the
formula above.  It should be noted that our convention here is the
opposite of that in \cite{ew3}).  There are canonical embeddings
$i_A:A\to M(A\rtimes_{\alpha,u}G)$ and $i_G:G\to
UM(A\rtimes_{\alpha,u}G)$ given on $f\in L^1(G,A)$ by
\begin{equation*}
\bigl(i_A(a)f\bigr)(t)=a\cdot f(t)\quad\text{and}\quad
\bigl(i_G(s)f\bigr)(t)=\alpha_s(f(s^{-1}t))u(s,s^{-1}t).
\end{equation*}
Then $(i_A,i_G)$ is a covariant homomorphism of $(A,G,\alpha,u)$ into
$M(A\rtimes_{\alpha,u}G)$, and for $f\in L^1(G,A)$, $i_A\times i_G(f)$
is the image of $f$ in $A\rtimes_{\alpha,u}G$ under the embedding of
$L^1(G,A)$ into its completion $A\rtimes_{\alpha,u}G$.

If $u\in Z^2\bigl(G, C(X,\T)\bigr)$, with $C(X,\TT)$ regarded as a
trivial $G$-module, then $(\id,u)$ is a twisted action of $G$ on
$C_{0}(X)$, and the central twisted transformation group algebras are
precisely the crossed products $C_0(X)\rtimes_{\id,u}G$ which arise in
this way.  For a good survey article on twisted transformation group
algebras we refer to \cite{judy-exp}.  Since we usually have
$\alpha=\id$ in this section, we shall often write $\hatu $ for the
dual action $\dualact{(\id, u)}$ of $\widehat{G}_{\ab}$.

To further reduce the notational overhead, we won't distinguish
between $v\in Z^2\bigl(G_{\ab}, C(X,\TT) \bigr)$ and its inflation,
$\inf v$ in $Z^2\bigl(G, C(X,\T)\bigr)$.  If $u\in Z^2\bigl(G,
C(X,\T)\bigr)$ is a product $u=v\cdot\sigma$ for some $v\in
Z^2_{\pt}(G_{\ab}, C(X,\TT))$ and $\sigma\in Z^2\bigl(G,
C(X,\T)\bigr)$ we want to obtain a description of
$C_0(X)\rtimes_{\id,u}G$ in terms of the algebras
$C_0(X)\rtimes_{\id,v}G_{\ab}$ and $C_0(X)\rtimes_{\id,\sigma}G$.  (Of
course this is motivated by Proposition \ref{prop-smooth}.)

We start by giving a description of $C_0(X)\rtimes_{\id,v}G_{\ab}$ in
terms of the bundle $p:\zmodule_v\to X$. The following result is well
known. It follows from the work of Smith \cite{hasm1} and the
discussion given in \cite[Remark 3.11]{90a}. However we feel it
worthwhile to include the following direct (and much shorter) proof.

\begin{lemma}\label{lem-pointwise}
   Suppose that $G$ is abelian and $v\in Z^2_{\pt}\bigl(G,
   C(X,\T)\bigr)$.  Then the dual system $(C_0(X)\rtimes_{\id,v}G,
   \widehat{G}, \hatv )$ is isomorphic to
   $(C_0(\zmodule_v),\widehat{G}, \tau)$, where $\tau:\widehat{G}\to
   \Aut C_0(\zmodule_v)$ is given by $\tau_{\gamma}
   (\psi)(f,x)=\psi(\overline{\gamma}\cdot f, x)$.
\end{lemma}
\begin{proof}
   Since $v$ is pointwise trivial, it is symmetric, i.e.,
   $v(s,t)=v(t,s)$ for every $s,t\in G$. Using this it is easy to check
   that convolution on $L^1\bigl(G, C(X)\bigr)$ is commutative.  Thus,
   $C_0(X)\rtimes_{\id,v}G$ is commutative, and isomorphic to
   $C_0\bigl((C_0(X)\rtimes_{\id,v}G)^{\wedge}\bigr)$ via the
   Gelfand-Naimark theorem.

   Thus we have to show that $\zmodule_v$ is
   $\widehat{G}$-equivariantly homeomorphic to
   $(C_0(X)\rtimes_{\id,u}G)^{\wedge}$. It is straightforward to check
   that the one-dimensional covariant representations are precisely the
   pairs $(\eps_x,f)$ with $(f,x)\in \zmodule_v$, where
   $\eps_x:C_0(X)\to\C$ denotes evaluation at $x$.  Thus we get a
   canonical bijection $\Phi:\zmodule_v\to
   (C_0(X)\rtimes_{\id,u}G)^{\wedge}$ given by
       \begin{equation*}
       (f,x)\mapsto
       \eps_x\rtimes f.
\end{equation*}
Since the action of a character $\gamma\in \widehat{G}$ on a covariant
representations $(U,\pi)$ is given by $(\gamma\cdot U,\pi)$, it
follows that $\Phi$ is $\widehat{G}$-equivariant.

So it only remains to check that $\Phi$ is continuous and open.  For
continuity, let $(f_n, x_n)$ converge to $(f,x)$ in $\zmodule_u$.
Then $\eps_{x_n}\rtimes f_n(h)$ converges to $\eps_x\rtimes f(h)$ for
all $h\in C_c(G\times X)$ by Lebesgue's dominated convergence theorem.
Since $C_c(G\times X)$ is dense in $C_0(X)\rtimes_{\id,u}G$ this
implies that $\eps_{x_n}\rtimes f_n$ converges to $\eps_x\rtimes f$ in
the weak-$*$ topology.  This proves continuity.

To prove openness, we could appeal to some deep results of Olesen and
Raeburn such as \cite[Corollary~2.3]{doir}.  Instead, we give a more
elementary argument.  We suppose that $\eps_{x_{n}}\rtimes f_{n}$
converges to $\eps_{x}\rtimes f$ in
$(C_{0}(X)\rtimes_{\id,u}G)^{\wedge}$.  In view of
Proposition~\ref{prop-bijection}, it will suffice to show that
$\set{(\underline{f}_{n},x_{n})}$ converges to $(\underline{f},x)$ in
$\zbarmodule_{v}$.  Let $h\in C_{c}(G\times X)$ be such that
$\eps_{x}\rtimes f(h)=1$.  Then if $\psi\in C_{c}(X)$, we note that
\begin{equation*}
     \eps_{x_{n}}\rtimes
        f_{n}\bigl(i_{C_{0}(X)}(\psi)h\bigr) =
\psi(x_{n})\bigl(\eps_{x_{n}}\rtimes
       f_{n}(h)\bigr).
\end{equation*}
Since $\eps_{x_{n}}\rtimes f_{n}\bigl(i_{C_{0}(X)}(\psi)h\bigr)$
converges to $\eps_{x}\rtimes f \bigl(i_{C_{0}(X)}(\psi)h\bigr)$, we
must have $\psi(x_{n})\to \psi(x)$ for all $\psi\in C_{c}(X)$.  Thus
$x_{n}\to x$.  Thus we may assume there is a $\psi\in C_{c}(X)$ such
that $\psi(x_{n})=1$ for all $n$. If $\varphi\in L^{1}(G)$ and we define $h\in
L^{1}\bigl( G, C_{0}(X) \bigr)$ by $h(s)=\varphi(s)\psi$, then for all
$\varphi\in L^{1}(G)$ we have
        \begin{equation}\label{eq:1}
          \int_{G}\varphi(s)f_{n}(s)\, d\mu(s)=\eps_{x_{n}}\rtimes f_{n}(h)
          \to \eps_{x}\rtimes
          f(h)=\int_{G}\varphi(s)f(s)\,d\mu(s).
        \end{equation}
        We claim \eqref{eq:1} implies that $\underline{f}_{n}\to
        \underline{f}$ in $\uC^{1}(G,\T)$.  In view of
        \cite[Proposition~6]{moore3}, it will suffice to show that
        \begin{equation}\label{eq:2}
        \int_{K}|f_{n}(s)-f(s)|\,d\mu(s)\to 0
        \end{equation}
        for each compact set $K\subset G$.  Since $|f(s)f_{n}(s)|=1$,
        we have
\begin{equation*}
    |f_{n}(s)-f(s)|^{2}=|1-\overline{f(s)}f_{n}(s)|^{2}\le 2 -
     2\operatorname{Re}\bigl(\overline{f(s)}f_{n}(s)\bigr).
\end{equation*}
But then by H\"older's inequality
\begin{align}
   \Bigl(\int_{K}|f_{n}(s)-f(s)|\,d\mu(s)\Bigr)^{2} &\le \mu(K)
   \int_{K} |f_{n}(s)-f(s)|^{2} \,d\mu(s)\notag \\
   &\le 2\mu(K)\operatorname{Re} \Bigl(
   \int_{K}1 -\overline{f(s)}f_{n}(s) \,d\mu(s)\Bigr)\notag \\
   &= 2\mu(K)\operatorname{Re} \Bigl(
   \int_{G}\charfcn_{K}(s)-\charfcn_{K}(s)\overline {f(s)}
   f_{n}(s)\,d\mu(s)\Bigr). \label{eq:4}
\end{align}
Since $\charfcn_{K}\cdot \overline{f}\in L^{1}(G)$, \eqref{eq:1}
implies that \eqref{eq:4} goes to zero, and the result follows.
\end{proof}

To state our main result in this section, we want to recall some
constructions from \cite[\S3]{ew3}.  A \emph{$C_0(X)$-algebra} is a
$C^*$-algebra $A$ equipped with a fixed nondegenerate $*$-homomorphism
$\Phi$ from $C_0(X)$ into the center $ZM(A)$ of the multiplier algebra
$M(A)$ of $A$.  This allows us to view $A$ as a $C_0(X)$-module, and
we shall usually write $\varphi\cdot a$ in place of $\Phi(\varphi)a$.
A twisted action $(\alpha,u)$ of $G$ on $A$ is called
\emph{$C_0(X)$-linear} if $\alpha_{s}(f\cdot a)=f\cdot \alpha_{s}(a)$
for all $s \in G$, $f\in C_0(X)$ and $a\in A$.  A $C_0(X)$-algebra $A$
should be thought of as an algebra of (semi-continuous) sections of a
bundle of $C^*$-algebras $A_x$, $x\in X$, where $A_x\cong
A/(C_0(X\setminus\set{x})\cdot A)$: the image of $a\in A$ in $A_{x}$
is denoted by $a(x)$.  Then $C_0(X)$-linearity means that $\alpha$
induces compatible twisted actions $(\alpha^x,u^x)$ on the fibres
$A_x$. The crossed product $A\rtimes_{\alpha,u}G$ of a $C_0(X)$-linear
twisted action is a $C_0(X)$-algebra with respect to the composition
$i_A\circ \Phi:C_0(X)\to M(A\rtimes_{\alpha,u}G)$, and the fibres are
just the crossed products $A_x\rtimes_{\alpha^x,u^x}G$.  In
particular, the dual action of $\widehat{G}_{\ab}$ on
$A\rtimes_{\alpha,u}G$ is again $C_0(X)$-linear and restricts to the
respective dual actions on the fibres $A_x\rtimes_{\alpha^x,u^x}G$.
It is worth mentioning that the dual $\ahat$ (respectively, the
primitive ideal space $\Prim A $) of a $C_0(X)$-algebra $A$ has an
induced bundle structure $q:\ahat \to X$ (resp. $q:\Prim A \to X$)
with fibres $\ahat _x$ (resp.  $\Prim A_x$).

Suppose that $L$ is an \emph{abelian} group, that $\alpha:L\to \Aut A
$ is a $C_0(X)$-linear (untwisted) action and that $p:Z\to X$ is a
locally compact free and proper $L$-bundle.  Since $L$ is abelian,
$\alpha^{-1}$ is an action and we can form the \emph{induced algebra}
$\Ind_L^Z(A,\alpha^{-1})$ which is the set of bounded continuous
functions $F:Z\to A$ satisfying
\begin{equation*}
\alpha_l(F(z))=F(l^{-1}\cdot z), \quad\text{for all $s\in G$ and
$z\in Z$, and}
\end{equation*}
such that $z\mapsto \|F(z)\|$ vanishes at infinity on $X=Z/L$.  As in
\cite{ew3}, we'll write $Z\times_{L}A$ in place of
$\Ind_L^Z(A,\alpha^{-1})$ to stress the analogy with classical
topological bundle constructions (see \cite[Definition~3.1(s)]{ew3}).
Note that $Z\times_LA$ is a $C^*$-algebra when equipped with the
pointwise operations and the supremum norm.  Moreover, $Z\times_LA$
carries a strongly continuous action $\Ind\alpha$ of $L$ given by
\begin{equation*}
(\Ind\alpha)_l(F)(z)=F(l^{-1}\cdot z)=\alpha_l(F(z)).
\end{equation*}

If $(A,L,\alpha)$ is a $C_0(X)$-system, then there is a $C_0(X\times
X)$-action on $Z\times_LA$ given by
\begin{equation*}
(h\cdot F)(z)(x)=h(p(z),x)F(z)(x), \quad h\in C_0(X\times X),
\end{equation*}
and the \emph{$L$-fibre product} $Z*A$ is defined as the restriction
of $Z\times_LA$ to the diagonal $\Delta(X)=\set{(x,x): x\in X}$.
Identifying $X$ with $\Delta(X)$ gives $Z*A$ the structure of a
$C_0(X)$-algebra and $\Ind\alpha$ restricts to a $C_0(X)$-linear
action $Z*\alpha$ of $L$ on $Z*A$.  Further details on this
construction and those in the previous two paragraphs can be found in
\cite[\S3]{ew3}.

Note that the construction of $Z*A$ is the $C^*$-algebraic analogue of
the usual construction of $L$-fibre products of topological
$L$-bundles given by
\begin{equation*}
Z*Y=\set{(z,y)\in Z\times Y: p(z)=q(y)}/L,
\end{equation*}
where $q:Y\to X$ is assumed to be a topological bundle over $X$
equipped with a compatible $L$-action on the fibres $Y_x$, and where
the quotient space is taken by the anti-diagonal action
$l\cdot(z,y)=(l z, l^{-1} y)$.  In particular, we always have
$(Z*A)^{\wedge}\cong Z*\ahat $ and $\Prim(Z*A)\cong Z*\Prim A $ with
respect to the bundle structures of $\ahat $ and $\Prim A $ (see
\cite[Proposition 3.5]{ew3}).

If $A$ and $B$ are $C_0(X)$-algebras, and if $\alpha:L\to \Aut A $ and
$ \beta:L\to \Aut B $ are $C_0(X)$-linear actions, then $(A,L,\alpha)$
and $(B, L,\beta)$ are \emph{$C_0(X)$-Morita equivalent} if there
exists an $A\sme B$ imprimitivity bimodule $\X$
satisfying\footnote{The left and right actions of $C_0(X)$ on $\X$ are
   obtained from extending the left and right actions of $A$ and $B$ to
   $M(A)$ and $M(B)$, respectively.}  $\varphi\cdot \xi=\xi\cdot
\varphi$ for all $\varphi\in C_0(X)$ and $\xi\in \X$, and such that
$\X$ carries a $C_0(X)$-linear action $\delta:L\to \Aut \X$ such that
\begin{equation*}
\blip A< \delta_l(\xi), \delta_l(\eta)>=
\alpha_l\bigl(\lip A<\xi,\eta>\bigr)\quad\text{and}
\quad \brip B< \delta_l(\xi), \delta_l(\eta)>=
\beta_l\bigl(\rip B<\xi,\eta>\bigr)
\end{equation*}
for all $\xi,\eta\in \X$ and $l\in L$.  Note that $C_0(X)$-Morita
equivalence implies equivalence of the topological $L$-bundles $\ahat
$ and $\bhat $ (resp. $\Prim A $ and $\Prim B $) --- see
\cite[Proposition~5.7]{rw-book}.

We can now state the main result of this section:

\begin{thm}[{cf., \cite[Theorem 5.3]{ew3}}]\label{thm-bundle}
   Suppose that $u$ has the form $v\cdot \sigma$ with $v\in Z^2_{\pt}(G_{\ab},
   C(X,\TT))$ and $\sigma\in Z^2\bigl(G, C(X,\T)\bigr)$.  Let $p:Z_v\to
   X$ be the free and proper $\widehat{G}_{\ab}$-bundle associated to
   $v$ as in Definition~\ref{def-Zu}. Then the systems
       \begin{equation*}
\bigl(C_0(X)\rtimes_{\id,u}G, \widehat{G}_{\ab},
       \hatu \bigr)\quad\text{and}\quad
       \bigl(Z_v*(C_0(X)\rtimes_{\id,\sigma}G), \widehat{G}_{\ab},
       \zmodule_v*\hats\bigr)
\end{equation*} are $C_0(X)$-Morita equivalent
systems.
\end{thm}

\begin{remark}\label{rem-iso}
   If $v$ is actually locally trivial, then a stronger result holds.
   It follows from \cite[Theorem 5.3]{ew3} that there exists a
   $\widehat{G}_{\ab}$-equivariant and $C_0(X)$-linear
   \emph{isomorphism} between $C_0(X)\rtimes_{\id,u}G$ and
   $Z_v*(C_0(X)\rtimes_{\id,\sigma}G)$.  The action $\delta$ appearing
   in that theorem, and which is used to compare $(\id, u)$ with $(\id,
   \sigma)$, is the locally unitary action $\delta:G\to \Aut C_0(X,\K)$
   corresponding to $\inf v\in Z^2_{\loc}\bigl(G, C(X,\T)\bigr)$ as
   constructed in \cite[Proposition 3.1]{HORR}.  However, the proof of
   \cite[Theorem 5.3]{ew3} requires localizations of both systems and
   breaks down if $v$ is only assumed to be pointwise trivial.  We do
   not know whether the stronger result of $C_0(X)$-linear isomorphism
   also holds in the situation of Theorem~\ref{thm-bundle} above.
\end{remark}

In what follows, we denote by $A\otimes_{C_0(X)}B$ the \emph{maximal}
balanced tensor product of $C_0(X)$-algebras $A$ and $B$ (see
\cite{blanchard2} and \cite[\S2]{ew3}).  It is obtained by restriction
of the $C_0(X\times X)$-algebra $A\otimes_{\max}B$ to the diagonal
$\Delta(X)$ and therefore carries a canonical structure as
$C_0(X)$-algebra.  If $(\alpha,u)$ and $(\beta,v)$ are $C_0(X)$-linear
twisted actions on $A$ and $B$, respectively, then the diagonal
twisted action $(\alpha\otimes\beta, u\otimes v)$ on
$A\otimes_{\max}B$ restricts to a $C_0(X)$-linear twisted action
$(\alpha\otimes_X\beta, u\otimes_Xv)$ on $A\otimes_{C_0(X)}B$ (see
\cite[\S4]{ew3} for more details).

Suppose now that $\alpha:L\to \Aut A $ is a $C_0(X)$-linear
(untwisted) action of the abelian group $L$, and let $p:Z\to X$ be a
free and proper $L$-bundle over $X$.  Let $\tau$ denote the action of
$L$ on $C_0(Z)$ given by $\tau_l(\psi)(z)=\psi(l^{-1}\cdot z)$.  Then
it follows that the crossed product $(C_0(Z)\otimes
A)\rtimes_{\tau\otimes\alpha^{-1}}L$ is a $C_0(X\times X)$-algebra,
and the restriction to the diagonal $\Delta(X)$ is isomorphic to
$(C_0(Z)\otimes_{C_0(X)}A)\rtimes_{\tau\otimes_X\alpha^{-1}}L$.  We
define a $C_0(X\times X)$-linear action $\underline{\alpha}$ of $L$ on
$(C_0(Z)\otimes A)\rtimes_{\tau\otimes\alpha^{-1}}L$ by the formula
\begin{equation*}
\underline{\alpha}_{\,l}(f)(s)=\id\otimes\alpha_l(f(s))
\quad\text{for  $f\in L^1(L, C_0(Z)\otimes_{C_0(X)}A)$ and $ l,s\in
    L$.}
\end{equation*}
Note that this extends to an automorphism of the crossed product since
$\id\otimes_X\alpha_l$ commutes with $\tau_h\otimes\alpha_{h^{-1}}$
for all $l,h\in L$. Since $\underline{\alpha}$ is $C_0(X\times
X)$-linear, it restricts to a $C_0(X)$-linear action
$\underline{\alpha}^X$ of $L$ on
$(C_0(Z)\otimes_{C_0(X)}A)\rtimes_{\tau\otimes_X\alpha^{-1}}L$.  The
proof of Theorem~\ref{thm-bundle} depends heavily on the following
result.

\begin{prop}\label{prop-equivalence}
   In the situation above, the systems
     \begin{equation*}
\bigl((C_0(Z)\otimes A)\rtimes_{\tau\otimes\alpha^{-1}}L,
     L,\underline{\alpha}\bigr)\quad\text{and}\quad \bigl(Z\times_LA, L,
     \Ind\alpha\bigr)
\end{equation*} are $C_0(X\times X)$-Morita equivalent, and the
systems
     \begin{equation*}
\bigl((C_0(Z)\otimes_{C_0(X)}A)\rtimes_{\tau\otimes_X\alpha^{-1}}L,
     L,\underline{\alpha}^X\bigr)\quad\text{and}\quad
     \bigl(Z*A, L, Z*\alpha\bigr)
\end{equation*}
are $C_0(X)$-Morita equivalent.
\end{prop}
\begin{proof} Note that it is enough to prove the first Morita
   equivalence, since the second will follow from the first by
   restricting to the diagonal $\Delta(X)$.

   The proof of the first equivalence is based on the proof of
   \cite[Theorem 2.2]{rw}.  Note that our algebra $Z\times_LA$ is equal
   to the algebra $GC(Z,A)^{\tau\otimes \alpha^{-1}}$ in the language
   of \cite{rw}.  The proof of \cite[Theorem 2.2]{rw} shows that we
   obtain a $(C_0(Z)\otimes A)\rtimes_{\tau\otimes\alpha^{-1}}L
   -Z\times_LA$ imprimitivity bimodule $\X$ by taking the completion of
   $C_c(Z,A)$ with respect to the left and right $(C_0(Z)\otimes
   A)\rtimes L$- and $Z\times_LA$-valued inner products and left and
   right actions of $C_c(L, C_0(Z)\otimes A)\subset (C_0(Z)\otimes
   A)\rtimes L$ and $Z\times_LA$ on $\X$ given by the formulas
\begin{align*}
   \lip C_c(L,C_0(Z)\otimes A)<\xi,\eta>(l,z)&=
   \xi(z)\alpha_{l^{-1}}\bigl(\eta(l^{-1}z)^*\bigr)\\
   \rip Z\times_LA< \xi,\eta>&=
   \int_L\alpha_{l}\bigl(\xi(lz)^*\eta(lz)\bigr)\,d\mu(l)\\
   f\cdot \xi(z)&=\int_Lf(l,z)\alpha_{l}\bigl(\xi(l z)\bigr)\,d\mu(l)\\
   \xi\cdot F(z)&=\xi(z)F(z),
\end{align*}
where $F\in Z\times_LA$, $f\in C_c(L, C_0(Z)\otimes A)$ and
$\xi,\eta\in C_c(Z,A)$. To see that these formulas are equivalent to
those given in \cite{rw}, note that our action $\alpha^{-1}$ plays the
role of the action $\beta$ in \cite{rw}, and that we may replace $l$
by $l^{-1}$ in all integrals above since $L$ is abelian and therefore
unimodular.

It is now easy to check that the $C_0(X\times X)$-actions on both
sides of $\X$ are given by the formula
\begin{equation*}
\varphi\cdot \xi(z)(x)=\varphi(p(z),x)\xi(z)(x),\quad \text{$\varphi\in
C_0(X\times X)$ and $ \xi\in C_c(Z,A)$.}
\end{equation*}
Moreover, if we define
$$\delta_l(\xi)(z)=\alpha_l(\xi(z)),\quad \text{for $\xi\in C_c(Z,A)$,
   $l\in L$ and $z\in Z$},$$
then it is straightforward to check that
$\delta$ extends to a $C_0(X\times X)$-linear action on the completion
$\X$ of $C_c(Z,A)$ which is compatible with the actions
$\underline{\alpha}$ and $\Ind\alpha$.
\end{proof}

Suppose that $(\alpha,u)$ is a twisted action of $G$ on $A$ and that
$N$ is a closed normal subgroup of $G$.  Then, depending on a choice
of a Borel section $c:G/N\to G$, Packer and Raeburn \cite[Theorem
4.1]{para1} showed that there is a twisted action $(\beta,w)$ of $G/N$
on the crossed product $A\rtimes_{\alpha,u}N$ such that
\begin{equation*}
A\rtimes_{\alpha,u}G\cong
(A\rtimes_{\alpha,u}N)\rtimes_{\beta,w}G/N.
\end{equation*}
We want to apply their result in the very special case where $G$ is
the direct product $L\times N$. In this case, we get a particularly
nice description of the twisted action $(\beta,w)$ and the isomorphism
$\Phi$:

\begin{prop}\label{prop-decom}
   Suppose that $G=L\times N$ is the direct product of two second
   countable locally compact groups $L$ and $N$ and let $(\alpha,u)$ be
   a twisted action of $G$ on $A$.  Then there is a twisted action
   $(\beta, w)$ of $L$ on $A\rtimes_{\alpha,u}N$ given by the formulas
     \begin{gather*}
       \beta_l(f)(n)=\alpha_l\bigl(f(n)\bigr)
       u\bigl((l,e),(e,n)\bigr)u\bigl((e,n),(l,e)\bigr)^*
       \quad \text{$f\in L^1(N,A)$, and} \\
       w(l,h)=i_A\bigl(u\bigl((l,e),(h,e)\bigr)\bigr),\quad l,h\in L.
     \end{gather*}
     With this action, there is an isomorphism between
     $A\rtimes_{\alpha,u}G$ and
     $\bigl(A\rtimes_{\alpha,u}N)\rtimes_{\beta,w}L$ which restricts to
     a \hm{} of $L^1$-algebras $\Phi:L^1(L\times N, A)\to L^1\bigl(L,
     L^1(N,A)\bigr)$ given by the formula
     \begin{equation*}
     \Phi(f)(l)(n)=
     f(l,n)u\bigl((e,n),(l,e)\bigr)^*.
\end{equation*}
\end{prop}
\begin{proof} The proof is basically a consequence of
   \cite[Theorem~4.1]{para1} --- in particular, the formula for the
   action $(\beta, w)$ directly follows from the formulas as given in
   \cite{para1} with respect to the cross section $L\to L\times N$
   defined by $ l\mapsto (l,e)$.  We only have to check that the
   isomorphism is given on the $L^1$-algebras by the above formula.
   For this let $(i_A, i_N)$ denote the canonical embeddings of $(A,G)$
   into $M(A\rtimes_{\alpha,u}N)$ and, similarly, let $(j_{A\rtimes N},
   j_L)$ denote the embeddings of $(A\rtimes_{\alpha,u} N, L)$ into
   $M\bigl((A\rtimes_{\alpha,u}N)\rtimes_{\beta,w}L\bigr)$.  Then it is
   shown on \cite[p.~307]{para1} that the pair $(k_A, k_G)$ defined by
\begin{equation*}
k_A=j_{A\rtimes N}\circ i_A\quad\text{and}\quad
k_G\bigl((l,n)\bigr)=j_{A\rtimes N}\bigl(i_A\bigl(u\bigl((e,n),
(l,e)\bigr)^*\bigr)
i_N(n)j_L(l)
\end{equation*}
is a covariant homomorphism of $(A,G,\alpha,u)$ into
$M\bigl((A\rtimes_{\alpha,u}N)\rtimes_{\beta,w}L\bigr)$ such that the
integrated form $k_A\rtimes k_G$ is the desired isomorphism.  Thus for
$f\in L^1(L\times N, A)$ we get
\begin{align*}
   k_A\rtimes k_G(f)&=\int_G k_A\bigl( f(l,n)\bigr)
   k_G\bigl((l,n)\bigr)\,d\mu_{L}(l,n)\\
   &=\int_L\int_N j_{A\rtimes N}\bigl(i_A\bigl(f\bigl((l,n)\bigr)
   u\bigl((e,n), (l,e)\bigr)^*\bigr)i_N(n)\bigr)
   j_L(l)\,d\mu_{N}(n)\, d\mu_{L}(l)\\
   &=\int_Lj_{A\rtimes N}\bigl(\Phi(f)(l)\bigr)j_L(l)\,d\mu_{L}(l)
   =j_{A\rtimes N}\rtimes j_L\bigl(\Phi(f)\bigr).
\end{align*}
This completes the proof.
\end{proof}

\begin{remark}\label{rem-equiv} Note that it follows directly
   from the formula of the isomorphism $A\rtimes_{\alpha,u}G\cong
   (A\rtimes_{\alpha,u}N)\rtimes_{\beta,w}L$ given in the proposition
   that this isomorphism is $\widehat{L}_{\ab}\times
   \widehat{N}_{\ab}$-equivariant. Moreover, if $A$ is a
   $C_0(X)$-algebra and $(\alpha,u)$ is $C_0(X)$-linear, then it
   follows from the definition of $(\beta, w)$, that it is
   $C_0(X)$-linear and the formula for the above isomorphism shows that
   it is $C_0(X)$-linear, too.
\end{remark}

We shall also need

\begin{prop}\label{prop-decom1}
   Suppose that $(\alpha,u)$ and $(\beta, v)$ are $C_0(X)$-linear
   twisted actions of $L$ and $N$ on the $C_0(X)$-algebras $A$ and $B$,
   respectively.  Define the twisted action $(\alpha\otimes_X\beta,
   u\otimes_Xv)$ of $L\times N$ on $A\otimes_{C_0(X)}B$ in the obvious
   way.  Then there exists a $C_0(X)$-linear and
   $\widehat{L}_{\ab}\times\widehat{N}_{\ab}$-equivariant isomorphism
   between
       \begin{equation*}
(A\otimes_{C_0(X)}B)\rtimes_{\alpha\otimes_X\beta,
       u\otimes_Xv}L\times N\quad\text{and}\quad
       \bigl(A\rtimes_{\alpha,u}L\bigr)\otimes_{C_0(X)}
       \bigl(B\rtimes_{\beta,v}N\bigr).
\end{equation*}
\end{prop}
\begin{proof} Restricting the twisted action
   $(\alpha\otimes_X\beta, u\otimes_Xv)$ to the subgroup $N\cong
   \set{e}\times N$ of $L\times N$ gives the action
   $(\id\otimes_X\beta, 1\otimes_Xv)$ of $N$ on $A\otimes_{C_0(X)}B$.
   It follows from \cite[Proposition 4.3]{ew3} that
   $(A\otimes_{C_0(X)}B)\rtimes_{\id\otimes_X\beta, 1\otimes_Xv}N$ is
   $C_0(X)$-linearly and $\widehat{N}_{\ab}$-equivariantly isomorphic
   to $A\otimes_{C_0(X)}(B\rtimes_{\beta,v}N)$.  Using the formula for
   this isomorphism as given \cite[Propsition 4.3]{ew3}, it follows
   from Proposition~\ref{prop-decom} that the decomposition action of
   $L$ on $(A\otimes_{C_0(X)}B)\rtimes_{\id\otimes_X\beta,
     1\otimes_Xv}N$ corresponds to the twisted action
   $(\alpha\otimes_X\id, u\otimes_X1)$ of $L$ on
   $A\otimes_{C_0(X)}(B\rtimes_{\beta,v}N)$.  The result then follows
   from another application of \cite[Proposition 4.3]{ew3}.
\end{proof}

We are now ready for the proof of Theorem~\ref{thm-bundle}.  The proof
relies on the above decomposition results, and the Takesaki-Takai
duality for twisted actions of abelian groups.

\begin{proof}[Proof of Theorem~\ref{thm-bundle}]
   We consider the diagonal twisted action $(\id\otimes_X\id,
   v\otimes_X\sigma)$ of $G_{\ab}\times G$ on
   $C_0(X)\otimes_{C_0(X)}C_0(X)\cong C_0(X)$.  If we restrict this
   action to the diagonal $\Delta(G)=\set{(\dot s,s):s\in G}\subset
   G_{\ab}\times G$ and identify $G$ with $\Delta(G)$ via $s\mapsto
   (\dot{s}, s)$, then it follows that the isomorphism
   $C_0(X)\otimes_{C_0(X)}C_0(X)\to C_0(X)$, given on elementary
   tensors by $\varphi\otimes\psi\to \varphi\cdot\psi$, carries
   $(\id_{\widehat{G}_{\ab}}\otimes_X\id_{G}, v\otimes_X\sigma)$ to the
   twisted action $(\id_{G}, v\cdot\sigma)=(\id_{G}, u)$.  Thus we get
   a natural $C_0(X)$-linear isomorphism
       \begin{equation*}
       C_0(X)\rtimes_{\id,u}G\cong
       \bigl(C_0(X)\otimes_{C_0(X)}C_0(X)\bigr)
       \rtimes_{\id\otimes_X\id, v\otimes_X\sigma}\Delta(G),
\end{equation*}
which transforms the dual action of $\widehat{G}_{\ab}$ to the dual
action of $\Delta(G)_{\ab}^{\wedge}$.

For the crossed product by the full group $G_{\ab}\times G$, it
follows from Proposition~\ref{prop-decom} that we have a
$C_0(X)$-linear and
$\widehat{G}_{\ab}\times\widehat{G}_{\ab}$-equivariant isomorphism
\begin{multline*}
   \bigl(C_0(X)\otimes_{C_0(X)}C_0(X)\bigr)
   \rtimes_{\id\otimes_{X}\id,v\otimes_X\sigma}G_{\ab}\times G
   \cong \\
   \bigl(C_0(X)\rtimes_{\id,v}G_{\ab}\bigr)\otimes_{C_0(X)}
   \bigl(C_0(X)\rtimes_{\id,\sigma}G\bigr).
\end{multline*}
By Lemma~\ref{lem-pointwise}, the algebra
$C_0(X)\rtimes_{\id,v}G_{\ab}$ is $\widehat{G}_{\ab}$-equivariantly
isomorphic to $C_0(\zmodule_v)$, and this isomorphism is clearly
$C_0(X)$-linear.  Thus we obtain a $C_0(X)$-linear and
$\widehat{G}_{\ab}\times\widehat{G}_{\ab}$-equivariant isomorphism
between
\begin{multline*}
   \bigl(C_0(X)\otimes_{C_0(X)}C_0(X)\bigr)
   \rtimes_{\id\otimes_{X}\id,v\otimes_X\sigma}G_{\ab}\times G\quad
   \text{and}\\
   C_0(\zmodule_v)\otimes_{C_0(X)}\bigl(C_0(X)\rtimes_{\id,\sigma}G).
\end{multline*}
We now split $G_{\ab}\times G$ as the product $G_{\ab}\times\Delta(G)$
via the isomorphism $(\dot{s}, t)\mapsto \bigl(\dot{t}^{-1}\dot{s},
(\dot{t},t)\bigr)$. Iterating the crossed product with respect to this
decomposition of $G_{\ab}\times G$ now provides $C_0(X)$-linear
isomorphisms
\begin{align*}
   C_0(\zmodule_v)\otimes_{C_0(X)}&\bigl(C_0(X)\rtimes_{\id,\sigma}G)\\
   &\cong\bigl(C_0(X)\otimes_{C_0(X)}C_0(X)\bigr)
   \rtimes_{\id\otimes_{X}\id,v\otimes_X\sigma}G_{\ab}\times G\\
   &\cong \bigl(\bigl(C_0(X)\otimes_{C_0(X)}C_0(X)\bigr)
   \rtimes_{\id\otimes_{X}\id,v\otimes_X\sigma}\Delta(G)\bigr)
   \rtimes_{\beta,w}G_{\ab}\\
   &\cong \bigl(C_0(X)\rtimes_{u}G\bigr)\rtimes_{\beta,w}G_{\ab}.
\end{align*}
We need to compare the natural
$\widehat{G}_{\ab}\times\widehat{G}_{\ab}$-action on
$\bigl(C_0(X)\rtimes_{u}G\bigr)\rtimes_{\beta,w}G_{\ab}$ with the
$\widehat{G}_{\ab}\times\widehat{G}_{\ab}$-action on
$C_0(\zmodule_v)\otimes_{C_0(X)}\bigl(C_0(X)\rtimes_{\id,\sigma}G)$
under the above isomorphism. Indeed, if we identify
$\widehat{G}_{\ab}$ with $\Delta(G)_{\ab}^{\wedge}$ via
$\chi(\dot{t},\dot t)=\chi(\dot t)$ (as we do in the last isomorphism
above), we see that our given isomorphism $\psi:G_{\ab}\times G\to
G_{\ab}\times \Delta(G)$ induces the isomorphism
$\widehat{\psi}:\widehat{G}_{\ab}\times\Delta(G)_{\ab}^{\wedge}\to
\widehat{G}_{\ab}\times\widehat{G}_{\ab}$ given by
\begin{equation*}
\widehat{\psi}(\gamma,\chi)(\dot{s},\dot t)=
\gamma(\dot{t}^{-1}\dot{s})\chi(\dot{t}).
\end{equation*}
It follows from this that the dual action $\hatu $ of
$\widehat{G}_{\ab}$ on $C_0(X)\rtimes_{\id,u}G$ (and then extended to
$\bigl(C_0(X)\rtimes_{u}G\bigr)\rtimes_{\beta,w}G_{\ab}$) corresponds
to the action $\id\otimes_X\hats $ of $\widehat{G}_{\ab}$ on
$C_0(\zmodule_v)\otimes_{C_0(X)}\bigl(C_0(X)\rtimes_{\id,\sigma}G)$,
while the dual action $(\beta,w)^{\wedge}$ of $\widehat{G}_{\ab}$ on
$\bigl(C_0(X)\rtimes_{u}G\bigr)\rtimes_{\beta,w}G_{\ab}$ corresponds
to the action $\tau\otimes_X\hats ^{-1}$ of $G_{\ab}$ on
$C_0(\zmodule_v)\otimes_{C_0(X)}\bigl(C_0(X)\rtimes_{\id,\sigma}G)$.
Thus it follows from Proposition~\ref{prop-equivalence} that we get a
$C_0(X)$-Morita equivalence between the systems
\begin{multline*}
   \bigl(\bigl(\bigl(C_0(X)\rtimes_{u}G\bigr)
   \rtimes_{\beta,w}G_{\ab}\bigr)
   \rtimes_{\dualact{(\beta,w)}}\widehat{G}_{\ab}, \;\widehat{G}_{\ab},
   \;
   \underline{\hatu }\bigr)\\
   \cong\bigl(\bigl(C_0(\zmodule_v)
   \otimes_{C_0(X)}\bigl(C_0(X)\rtimes_{\id,\sigma}G)\bigr)
   \rtimes_{\tau\otimes\hats ^{-1}}\widehat{G}_{\ab},\;
   \widehat{G}_{\ab},\; \underline{\hats }^X\bigr)
\end{multline*}
and
\begin{equation*}
\bigl(\zmodule_v*\bigl(C_0(X)\rtimes_{\id,\sigma}G\bigr),\;\widehat{G}_{\ab},\;
\zmodule_v*\hats \bigr),
\end{equation*}
where $\underline{\hatu }$ denotes the canonical action induced by
$\hatu $ on
$\bigl(\bigl(C_0(X)\rtimes_{u}G\bigr)\rtimes_{\beta,w}G_{\ab}\bigr)
\rtimes_{\dualact{(\beta,w)}}\widehat{G}_{\ab}$.  Now the
Takesaki-Takai theorem for twisted actions (see
\cite[Theorem~3.1]{quigg86}) implies that
\begin{equation*}
\bigl(\bigl(C_0(X)\rtimes_{u}G\bigr)\rtimes_{\beta,w}G_{\ab}\bigr)
    \rtimes_{\dualact{(\beta,w)}}\widehat{G}_{\ab}
    \cong \bigl(C_0(X)\rtimes_{\id,u}G\bigr)\otimes\K(L^2(G_{\ab})),
\end{equation*}
and this isomorphism carries the action $\underline{\hatu }$ to $\hatu
\otimes\id_{\K}$. Since the systems
    \begin{equation*}
\bigl(\bigl(C_0(X)\rtimes_{\id,u}G\bigr)\otimes \K(L^2(G_{\ab})),\;
    \widehat{G}_{\ab},\;\hatu \otimes\id_{\K}\bigr)\quad\text{and}\quad
    \bigl(C_0(X)\rtimes_{\id,u}G,\;\widehat{G}_{\ab},\;\hatu \bigr)
\end{equation*}
are clearly $C_0(X)$-Morita equivalent, the result follows.
\end{proof}

We are now going to use Theorem~\ref{thm-bundle} to give a bundle
theoretic description of $C_0(X)\rtimes_{\id,u}G$ when $G$ is smooth.
Then Proposition~\ref{prop-smooth} implies that we obtain a
factorization $u=v\cdot u_{\varphi}$, where $v\in
Z^2_{\pt}\bigl(G_{\ab}, C(X,\TT)\bigr)$ and $u_{\varphi}$ is obtained
by pulling back a given cocycle $\eta\in Z^2(G,Z)$ corresponding to a
representation group $1\to Z\to H\to G\to 1$ for $G$ via the
continuous map $\varphi:X\to H^2(G,\TT)\cong \widehat{Z}$ defined by $
x\mapsto [u(x)]$.  Recall from \cite{rw, ew3} that if $A$ is a
$C_0(Y)$-algebra and if $\varphi:X\to Y$ is a continuous map, then the
pull-back $\varphi^*(A)$ is defined as the balanced tensor product
$C_0(X)\otimes_{C_0(Y)}A$, where $C_0(X)$ is viewed as a
$C_0(Y)$-algebra via $\varphi:X\to Y$.  Note that $\varphi^*(A)$
becomes a $C_0(X)$-algebra via the canonical embedding $C_0(X)\to
M\bigl(C_0(X)\otimes_{C_0(Y)}A\bigr)$.  Moreover, if $\alpha:L\to \Aut
A $ is a $C_0(Y)$-linear action on $A$, then
$\varphi^*(\alpha)=\id\otimes_Y\alpha$ is a $C_0(X)$-linear action on
$\varphi^*(A)$.  The following description of
$C_0(X)\rtimes_{\id,u_{\varphi}}G$ follows from \cite[Lemmas 6.3~and
6.5]{ew3}.

\begin{prop}\label{prop-pull-back}
   Let $1\to Z\to H\to G\to 1$ be a representation group for $G$. Let
   $C^*(H)$ be viewed as a $C_0(\widehat{Z})$-algebra via the canonical
   embedding $C_0(\widehat{Z})\cong C^*(Z)\to ZM(C^*(H))$ given by
   convolution.  Let $\varphi:X\to H^2(G,\TT)\cong \widehat{Z}$ be a
   continuous map, and let $u_{\varphi}\in Z^2\bigl(G, C(X,\T)\bigr)$
   be the cocycle defined in Proposition~\ref{prop-smooth} (with
   respect to \emph{any} cocycle $\eta\in Z^2(G,Z)$ corresponding to
   $H$).  Further, let $\delta$ denote the dual action of
   $\widehat{G}_{\ab}=\widehat{H}_{\ab}$ on $C^*(H)$.  Then the systems
     \begin{equation*}
\bigl(C_0(X)\rtimes_{\id, u_{\varphi}}G,\;\widehat{G}_{\ab},\;
     \hatu _{\varphi}\bigr)
     \quad\text{and}\quad
     \bigl(\varphi^*\bigl(C^*(H)\bigr),\;\widehat{G}_{\ab},\;
     \varphi^*(\delta)\bigr)
\end{equation*}
are $C_0(X)$-isomorphic.
\end{prop}

We can now gather our results to obtain a general description of the
bundle structure of $C_0(X)\rtimes_{\id,u}G$ in terms of a given
representation group $H$ for $G$.

\begin{thm}\label{thm-general}
   Suppose that $G$ is smooth with representation group $1\to Z\to H\to
   G\to 1$, and that $u\in Z^2\bigl(G, C(X,\T)\bigr)$. Let
   $u_{\varphi}\in Z^2\bigl(G, C(X,\T)\bigr)$ be as above with
   $\varphi(x)=[u(x)]$ for all $x\in X$.  If $v:= u\cdot
   \overline{u}_{\varphi}$, then $v\in Z^2_{\pt}(G, C(X,\TT))$ and
   there exists a $C_0(X)$-Morita equivalence between the systems
\begin{equation*}
\bigl(C_0(X)\rtimes_{\id,u}G,\;\widehat{G}_{\ab},\;\hatu \bigr)
\quad\text{and}\quad
\bigl(\zmodule_v*\varphi^*\bigl(C^*(H)\bigr),\;\widehat{G}_{\ab},
\;\zmodule_v*\bigl(\varphi^*(\delta)\bigr)\bigr).
\end{equation*}
\end{thm}
\begin{proof} The proof is now a direct consequence of
   Propositions \ref{prop-smooth}~and \ref{prop-pull-back}, and
   Theorem~\ref{thm-bundle}.
\end{proof}

\begin{remark}\label{rem-final}
   (1) If the cocycle $v=u\cdot\overline{u}_{\varphi}$ in the above
   theorem is actually \emph{locally trivial}, then it follows from
   Remark~\ref{rem-iso} that the $C_0(X)$-Morita equivalence in the
   theorem can be replaced by $C_0(X)$-isomorphism.  By
   \cite[Theorem~2.1]{ros2}, this is automatically the case whenever
   $\widehat{G}_{\ab}$ is compactly generated.

   (2) If $A$ is a $\CR(X)$-algebra in the sense of \cite{ew2,ew3}
   (e.g., if $A$ is unital and $X$ is the complete regularization of
   $\Prim A $ as in \cite[Definition 2.5]{ew2}), then any inner action
   of $G$ on $A$ determines a unique class $[u]\in H^2\bigl(G,
   C(X,\T)\bigr)$ (see \cite[\S 0]{rr} and \cite[\S2]{ew2}). It is
   shown in \cite[Corollary 4.7]{ew3} that the crossed product
   $A\rtimes_{\alpha}G$ is then $C_0(X)$-linearly and
   $\widehat{G}_{\ab}$-equivariantly isomorphic to
   $\bigl(C_0(X)\rtimes_{\id,u}G\bigr)\otimes_{C_0(X)}A$, so
   Theorem~\ref{thm-general} also gives new insights into the structure
   of crossed products by inner actions.
\end{remark}

We end this section by giving a description of
$C_0(X)\rtimes_{\id,u}G$ where $u$ is a cocycle in $Z^2\bigl(G,
C(X,\T)\bigr)$ with constant evaluation map $[u(x)]=[\omega]\in
H^2(G,\TT)$.  Note that if $H^2(G,\TT)$ is discrete (as considered by
Smith in \cite{hasm3}), then every cocycle $u\in Z^2\bigl(G,
C(X,\TT)\bigr)$ has a direct sum decomposition into cocycles $u_i\in
Z^2\bigl(G, C(X_i,\TT)\bigr)$ such that $X$ is a disjoint union of the
clopen subsets $X_i\subseteq X$, and each $u_i$ has constant evaluation
map. It is then easy to see that we get a decomposition
$$C_0(X)\rtimes_{\id,u} G\cong \bigoplus_i C_0(X_i)\rtimes_{\id,
  u_i}G.$$
As is standard, $C^*(G,\omega)$ will denote the twisted
group algebra $\CC\rtimes_{\id,\omega}G$ of $G$ with respect to
$\omega\in Z^2(G,\TT)$.  Of course, if $\omega$ is trivial, then
$C^*(G,\omega)$ is the full group $C^*$-algebra $C^*(G)$ of $G$.

\begin{thm}\label{thm-pt}
   Assume that $u\in Z^2\bigl(G, C(X,\TT)\bigr)$ has constant evaluation
map $x\mapsto [u(x)]:=[\omega]\in H^2(G,\TT)$. Suppose further that 
$v=u\cdot \overline{\omega}\in Z^2_{\pt}\bigl(G, C(X,\T)\bigr)$ 
satisfies one of
   the equivalent conditions of Proposition~\ref{prop-lifting}
(which is automatic if $H^2(G,\TT)$ is Hausdorff). Then
   the system
   $\bigl(C_0(X)\rtimes_{\id,u}G,\;\widehat{G}_{\ab},\;\hatv\bigr)$ is
   $C_0(X)$-Morita equivalent to
   $\bigl(\zmodule_v\times_{\widehat{G}_{\ab}}C^*(G,\omega),\;
   \widehat{G}_{\ab},\;\Ind\delta\bigr)$, where $\delta:G_{\ab}\to
   \Aut(C^*(G,\omega))$ denotes the dual action and the $C_0(X)$-structure of
   $\zmodule_v\times_{\widehat{G}_{\ab}}C^*(G,\omega)$ is given by
   $\bigl(\psi\cdot F\bigr)(z)=\psi(p(z))F(z)$.
\end{thm}
\begin{proof}
   If we apply Theorem~\ref{thm-bundle} to the factorization $u=v\cdot \omega$,
we obtain a $\widehat{G}_{\ab}$-equivariant $C_0(X)$-Morita
   equivalence between $C_0(X)\rtimes_{\id,u}G$ and
\begin{equation*}
\zmodule_v*\bigl(C_0(X)\rtimes_{\id,\omega}G\bigr)\cong
\zmodule_v*\bigl(C_0(X)\otimes C^*(G,\omega)\bigr)\cong
\zmodule_v\times_{\widehat{G}_{\ab}}C^*(G,\omega)
\end{equation*}
(with respect to the obvious identifications),
where the last isomorphism follows from \cite[Remark 3.4(c)]{ew3}.
\end{proof}

If $v\in Z^2_{\loc}\bigl(G, C(X,\T)\bigr)$, then the conditions of
Proposition~\ref{prop-lifting} are automatically satisfied
(Remark~\ref{rem-local}), and Remark~\ref{rem-iso} implies that for
such $u$ we may replace $C_0(X)$-Morita equivalence by
$C_0(X)$-isomorphism in the statement of Theorem~\ref{thm-pt}.

\section{The Group $C^*$-Algebras of Central Group
    Extensions}\label{sec:group}

In this section we use our methods to study of the group \cs-algebra
$C^*(L)$ of a central extension $1\to N\to L\to G\to 1$.  Of course,
the study of such algebras and their dual spaces is one of the main
motivations for studying central twisted transformation group
algebras.

To each central extension as above, we can associate a cocycle
$\eta\in Z^2(G,N)$ of the form $\eta=\partial c$ for a Borel
cross-section $c:G\to L$ satisfying $c(e)=e$.
Viewing $N$ as the dual of $\widehat N$, $\eta$ can be viewed as a
cocycle in $Z^2(G, C(\widehat{N},\TT))$, and it follows from
\cite[Theorem~4.1]{para1} (but see also
\cite[Lemma 6.3(a)]{ew3}), that $C^*(L)$ is $C_0(\widehat{N})$-linearly
and $\widehat{G}_{\ab}$-equivariantly isomorphic to
$C_0(\widehat{N})\rtimes_{\id,\eta}G$, where the
$C_0(\widehat{N})$-structure on $C^*(L)$ is given by the canonical
inclusion $C_0(\widehat{N})\cong C^*(N)\to ZM(C^*(L))$ given by
convolution (compare with the discussion preceding Proposition
\ref{prop-pull-back}). Thus the results of the preceding sections are
directly applicable to the study of $C^*(L)$. However, for central
twisted transformation group algebras associated to central group
extensions,
many of the abstract constructions in the preceding
sections, like the bundle $\zmodule_v$, can be realized
quite naturally on the group level (e.g., Remark~\ref{rem-pt-group}(5)).

\begin{definition}\label{def-pt-group}
    A central extension $1\to N\to L\to G\to 1$ of $G$ by $N$ is called
    \emph{pointwise trivial} if every character $\chi\in \widehat{N}$
    can be extended to a character of $L$; that is, if the restriction
    map $\res:H^1(L,\TT)\to H^1(N,\TT)=\widehat{N}$ is surjective.  We
    denote by $Z^2_{\pt}(G,N)$ the cocycles, and by $H^2_{\pt}(G,N)$ the
    classes in $H^2(G,N)$, corresponding to pointwise trivial extensions.
\end{definition}

\begin{remark}\label{rem-pt-group} We collect some
    straightforward observations on pointwise trivial extensions, which
    are no doubt well known to the experts.

    (1) A central extension $1\to N\to L\to G\to 1$ is pointwise trivial
    if and only if any corresponding cocycle $\eta$, viewed as an
    element of $Z^2(G, C(\widehat{N},\TT))$ is pointwise trivial in the
    sense of Definition~\ref{def-pointwise}.  This follows directly from
    the Hochschild-Serre exact sequence
       \begin{equation*}
\xymatrixnocompile@1{%
1\ar[r]& H^1(G,\TT)\ar[r]^{\inf} &
       H^1(L,\TT)\ar[r]^{\res} & H^1(N,\TT)
       \ar[r]^{\tg}  &H^2(G,\TT)}
\end{equation*}
(see \cite[Chap I \S5]{moore1}), but can easily be computed
directly.

(2) If $G$ is abelian, then the pointwise unitary extensions of $G$ by
$N$ are precisely the abelian extensions.  Indeed, if $1\to N\to L\to
G\to 1$ is a pointwise trivial extension with $G$ abelian, then one
easily checks that the characters of $L$ separate the points of $L$.
It follows that $[L,L]$ is trivial and $L$ is abelian.
Thus if $G$ is abelian, then
$H^2_{\pt}(G,N)=H^2_{\ab}(G,N)$,
where $H^2_{\ab}(G,N)$ denotes the set of cohomology classes
corresponding to the abelian extensions.

(3) If $1\to N\to L\to G\to 1$ is a pointwise trivial extension, then
the quotient map $L\mapsto L_{\ab}$ is injective on $N$ (since the
characters of $L$ separate the points of $N$). Thus we obtain an
abelian exact sequence $1\to N\to L_{\ab}\to G_{\ab}\to 1$, and the
extension $1\to N\to L\to G\to 1$ is actually inflated from this
abelian extension.  Recall that if $1\to N\to M\stackrel{p}{\to}
G_{\ab}\to 1$ is any extension of $G_{\ab}$, then the inflation of
this extension is the extension $1\to N\to \inf(M)\to G\to 1$ obtained
as follows.
We set $\inf(M)=\set{(m,s)\in M\times G:
p(m)=q(s)}$, where $q:G\to G_{\ab}$ is the quotient map.  The
inclusion $N\to \inf{M}$ sends $ n\mapsto (n,e)$,
and the quotient map $\inf{M}\to G$ sends $ (m,s)\mapsto s$.  The
isomorphism $L\cong\inf(L_{\ab})$ is given by $l\mapsto \bigl(p(l),
q(l)\bigr)$, where $p:L\to L_{\ab}$ and $q:L\to G$ are the quotient
maps.

(4) Of course, inflation of extensions in the above sense corresponds
to the inflation of the corresponding group cocycles. Indeed, if $1\to
N\to M\to G_{\ab}\to 1$ is as above and if
$c:G_{\ab}\to M$ is a Borel
section, then we obtain a Borel section $d:G\to \inf(M)$ by defining
$d(s)=\bigl(c(q(s)),s\bigr)$. Of course we then get $\partial_G
d=\inf\partial_{G_{\ab}} c$.  Combining this with~(3), we see that
inflation determines an isomorphism $\inf:H^2_{\ab}(G_{\ab}, N)\to
H^2_{\pt}(G,N)$.

(5) If $\eta\in Z^2_{\pt}(G, C(\widehat{N},\TT))$ is a cocycle
corresponding to a pointwise trivial extension $1\to N\to L\to G\to
1$, then the $\widehat{G}_{\ab}$-bundle $p:\zmodule_{\eta}\to
\widehat{N}$ is isomorphic to the bundle
$\res:\widehat{L}_{\ab}=H^1(L,\TT)\to\widehat{N}$.  Indeed, if
$\eta=\partial c$ for some cross section $c:G\to L$, then we define a
map $\Psi: \zmodule_{\eta}\to
\widehat{L}_{\ab}$ by
\begin{equation*}
\Psi(f,\chi)(c(s)n)=f(s)\chi(n).
\end{equation*}
This is well defined, since
\begin{align*}
    \partial_L(\Psi(f,\chi))(c(s)n, c(t)m)
    &=f(s)\chi(n)f(t)\chi(m)\overline{f(st)\chi(c(st)^{-1}c(s)c(t)nm)}\\
    &=\partial_G(f)(s,t)\overline{\eta(s,t)(\chi)}=1,
\end{align*}
so $\Psi(f,\chi)\in \widehat{L}_{\ab}$. Since pointwise convergence of
characters implies uniform convergence on compact sets (see
\cite[Theorem~8]{moore3}), the map $\Phi$
is continuous, and it is certainly $\widehat{G}_{\ab}$-equivariant.
The assertion follows from the fact that $\widehat{L}_{\ab}\to
\zmodule_{\eta}$ defined by $ \mu\mapsto (\mu\circ c, \mu|_N)$ is a
continuous inverse for $\Psi$.
\end{remark}

Thus, as a direct corollary of item~(5) of the above remark and
Theorem~\ref{thm-pt} we obtain

\begin{cor}\label{cor-pt-group}
    Suppose that $1\to N\to L\to G\to 1$ is a pointwise trivial central group
    extension and let $\delta:\widehat{G}_{\ab}\to \Aut C^*(G )$ denote
    the dual action.  Then $C^*(L)$ is
    $C_0(\widehat{N})$-Morita equivalent to
    $\widehat{L}_{\ab}\times_{\widehat{G}_{\ab}}C^*(G)\cong
    \Ind_{\widehat{G}_{\ab}}^{\widehat{L}_{\ab}}(C^*(G),\delta^{-1})$.
    (In fact, the corresponding $\widehat{G}_{\ab}$-systems are
    $C_0(\widehat{N})$-Morita equivalent.)
\end{cor}

Again, if the bundle $\res:\widehat{L}_{\ab}\to \widehat{N}$ is
locally trivial (which is automatic if $G_{\ab}$ is compactly
generated), then we can replace $C_0(\widehat{N})$-Morita equivalence
by $C_0(\widehat{N})$-isomorphism.

We are now going to discuss the group algebra of general central
extensions of
a smooth group $G$ with a fixed representation group $1\to Z\to H\to
G\to 1$. Let $\mu\in Z^2(G,Z)$ be a corresponding cocycle.  Then,
identifying $H^2(G,\TT)$ with $\widehat{Z}$, the transgression map for
$1\to N\to L\to G\to 1$ (which is just the evaluation map $\chi\mapsto
[\eta(\chi)]\in H^2(G,\TT)$, if $\eta\in Z^2(G,C(\widehat{N},\TT))$ is
a cocycle corresponding to the given extension) determines a
homomorphism $\varphi:\widehat{N}\to \widehat{Z}$.  Let
$\widehat{\varphi}:Z\to N$ denote the dual homomorphism defined by
$\widehat{\psi}(z)(\chi)=z(\psi(\chi))$ (where we identify $Z$ with
the dual of $\widehat{Z}$ and $N$ with the dual of $\widehat{N}$ via
Pontryagin duality).  Then we obtain a cocycle
$\widehat{\varphi}_*(\mu)\in Z^2(G,N)$ by defining
$\widehat{\varphi}_*(\mu)(s,t)=\widehat{\varphi}\bigl(\mu(s,t)\bigr)$.
A short computation shows that this cocycle, viewed as a cocycle in
$Z^2\bigl(G, C(\widehat{N},\TT)\bigr)$, is precisely the one we obtain
from the evaluation map for $\eta$ via the process described in
Proposition~\ref{prop-smooth}.  It follows in particular that
$\eta\cdot \widehat{\varphi}_*(\mu)^{-1}\in Z^2_{\pt}(G,N)$.  Thus, a
small variation on the proof of Proposition~\ref{prop-smooth} gives us

\begin{prop}\label{prop-smoothgr}
    Let $1\to Z\to H\to G\to 1$ be a representation group for $G$ and
    let $\mu\in Z^2(G,Z)$ be a corresponding cocycle.  Then, for any
    locally compact abelian group $N$, viewed as a trivial $G$-module, we
    get an isomorphism
      \begin{equation*}
H^2_{\ab}(G_\ab,N)\times \operatorname{Hom}(Z,N)\cong H^2(G,N)
\end{equation*}
defined by $ \bigl([\eta], \psi\bigr)
      \mapsto \bigl[\inf\eta\cdot \psi_*(\mu)\bigr]$.
\end{prop}

As a consequence of the above discussion and
Theorem~\ref{thm-general} we obtain

\begin{thm}\label{thm-groupex}
    Let $1\to N\to L\to G\to 1$, $\eta\in Z^2(G,N)$,
    $\varphi:\widehat{N}\to \widehat{Z}$ and $\eta\cdot
    \widehat{\varphi}_*(\mu)\in Z^2_{\pt}(G,N)$ be as in the discussion
    preceding Proposition~\ref{prop-smoothgr}, and let $1\to N\to\Lprime\to G\to
    1$ be the central extension corresponding to the pointwise trivial
    cocycle $\eta\cdot \widehat{\varphi}_*(\mu)^{-1}$.  Then $C^*(L)$ is
    $C_0(\widehat{N})$-Morita equivalent to
    $\widehat{\Lprime}_{\ab}*\varphi^*(C^*(H))$. (Again, if we consider
    the $\widehat{G}_{\ab}$-actions, the Morita equivalence passes to
    the dynamical systems.)
\end{thm}

It is actually easy to give a direct construction of the pointwise
trivial extension $1\to N\to \Lprime\to G\to 1$ corresponding to the
cocycle $\eta\cdot \widehat{\varphi}_*(\mu)^{-1}$ without even
mentioning the cocycles. For this let $1\to N\to
L\stackrel{p}{\to}G\to1$ be the original extension corresponding to
$\eta$.  Let $q:H\to G$ denote the quotient map for the representation
group $H$. Define
\begin{equation*}
\Lprime=\set{(l,h)\in L\times H: p(l)=q(h)}/\Delta(Z),
\end{equation*}
where $\Delta(Z)=\set{(\widehat{\varphi}(z), z):z\in Z}$. Then we obtain
a central extension
\begin{equation*}
\xymatrixnocompile@1{%
    1 \ar[r]& N \ar[rr]^{n\mapsto [n,e]}&&
     \Lprime \ar[rr]^{[l,h]\mapsto p(l)} && G \ar[r]& 1.}
\end{equation*}
We claim that this extension corresponds to the cocycle
$\eta\cdot\widehat{\varphi}_*(\eta)$ of the theorem. Indeed, if we
   choose Borel sections $c:G\to L$ and $d:G\to H$ such that
$\eta=\partial c$ and $\mu=\partial d$, then we get a Borel section
$c\times d:G\to \Lprime$ by defining $(c\times d)(s)=\bigl[c(s),d(s)\bigr]$.
We then compute
\begin{align*}
   \partial(c\times d)(s,t)&=\bigl[c(s),d(s)\bigr]\bigl[c(t),d(t)\bigr]
   \bigl[c(st),d(st)\bigr]^{-1}\\
   &=\bigl[\eta(s,t), \mu(s,t)\bigr]\\
   &=\bigl[\eta(s,t)\widehat{\varphi}(\mu(s,t))^{-1}, e\bigr],
\end{align*}
which clearly proves the claim.

We finish with some straightforward examples which illustrate our
results. 

\begin{example}\label{ex-zz2}
  Let $G=\ZZ^2$. Recall that the discrete Heisenberg group $H_d$ is
  the set $\ZZ^3$ with multiplication given by 
  $$(n_1, m_1,
  l_1)(n_2,m_2, l_2)=(n_1+n_2, m_1+m_2, l_1+l_2+n_1m_2).$$
  It is easy
  to check that
  \begin{equation*}
    \xymatrix{1\ar[r]& \protect{\ZZ}\ar[r] & H_{d}\ar[r]& \protect{\ZZ}^{2} \ar[r] & 1}
  \end{equation*}
     is a representation group for $G$
     and that $\widehat{\ZZ}=\TT\cong H^2(\ZZ^2,\TT)$ via
     $z\mapsto[\omega_z]$, where
     $$\omega_z\big((n_1,m_1),(n_2,m_2)\big)=z^{n_1m_2}.$$
     Since every abelian extension of $\ZZ^2$ by some
      group $N$ splits, it follows from
     Proposition~\ref{prop-smoothgr} that
     $H^2(\ZZ^2,N)\cong \operatorname{Hom}(\ZZ,N)=N$. Thus each $n\in
     N$ determines a central extension
     \begin{equation*}
       \xymatrix{1\ar[r]& N\ar[r] & L\ar[r]& \protect{\ZZ}^{2} \ar[r] & 1,}
     \end{equation*}
     and Theorem \ref{thm-groupex} gives an isomorphism
     between $C^*(L)$ and the pull-back  $n^*(C^*(H_d))$
     (where we identify $n$ with the character of $\widehat{N}$
     given by evaluation).
     Recall from \cite{andpas:hjm89} that $C^*(H_d)$ is a continuous bundle
     over $\TT$ with fibers $A_z$ given by the
     rotation algebras $A_{\theta}$ where
     $z=e^{2\pi i\theta}$ for some $\theta\in [0,1]$.
     Hence $C^*(L)$ is a continuous bundle over $\widehat{N}$
     with fibre $A_{\chi(n)}$ at the base point $\chi\in \widehat{N}$,
whose global structure is completely determined by the global
structure of $C^*(H_d)$ as a bundle over $\widehat{\ZZ}=\TT$.
\end{example}

In \cite{ew2} we gave explicit constructions for the representation
groups for $\ZZ^n$ and $\RR^n$. Using these, we can also apply our
results to 
these groups. In all these cases,
the abelian extensions vanish, so that we have
$H^2(G, N)\cong \operatorname{Hom}(Z, N)$, where $Z$ denotes the center of
the corresponding representation group $H$. Therefore, the group
algebras of central extensions of $\ZZ^n$ and $\RR^n$ by $N$ are isomorphic
to the pull-backs of $C^*(H)$ via the 
corresponding dual maps $\varphi:\widehat{N}\to
\widehat{Z}$. 

\begin{example}\label{ex-nontrivial}
Let $G=\ZZ\times\ZZ_2$, where $\ZZ_2=\ZZ/2\ZZ$. Then an easy 
application of \cite[Proposition 4.5]{ew3} shows that
$H=\ZZ\times \ZZ_2\times \ZZ_2$ with multiplication given by
$$(n, [i], [j])\cdot (m, [l],[k])= (n+m, [i+l], [j+k+n\cdot l]), \quad
n,m,i,j,k,l\in \ZZ,$$
is a representation group for $G$ with center
$Z=\ZZ_2$. In particular, we have
$H^2(G,\TT)\cong\widehat{\ZZ_2}=\ZZ_2$. As there are many non-trivial
abelian extensions of $G$ by locally compact abelian groups $N$, we
have, in general, a nontrivial decomposition
$H^2(G,N)=H^2_{\ab}(G,N)\oplus \operatorname{Hom}(Z,N)$. It is then an
straightforward exercise to apply our results to the group algebras of
the corresponding central extensions of $G$ by $N$.
\end{example}

It is also interesting to revisit Example \ref{ex-open}
to illustrate some differences between the general situation
of central twisted crossed products compared to central 
group extensions for possibly non-smooth groups.

\begin{example}\label{ex-final}
Let $G=\RR^2\times \TT^2$ with
multiplication given by
\begin{equation*}
(s_1, t_1, z_1, w_1)(s_2, t_2, z_2, w_2) =(s_1+s_2, t_1+t_2,
    e^{is_1t_2}z_1z_2, e^{i\theta s_1 t_2}w_1w_2),\end{equation*}
    where $\theta$ is any fixed irrational real number.
In Example \ref{ex-open} we constructed a pointwise
unitary cocycle
$u\in Z^2\bigl( G, C(X,\TT)\bigr)$, with $X=\{\frac{1}{n}:n\in
\NN\}\cup\{0\}$,  which is not inflated
from $G_{\ab}=\RR^2$. In particular,
$H^2_{\pt}\bigl(G, C(X,\TT)\bigr)\neq
H^2_{\pt}\bigl(\RR^2, C(X,\TT)\bigr)=\{0\}$.
On the other hand, it follows from part (4)
of Remark \ref{rem-pt-group} that for every abelian locally compact
group $N$ we do have an inflation isomorphism
$$
  \inf:   H^2_{\ab}(\RR^2, N) = H^2_{\ab}(G_{\ab},N)\to
H^2_{\pt}(G,N),
$$
from which it follows that $H^2_{\pt}(G,N)=\{0\}$ for all $N$.
Thus, the transgression map
$\tg:\widehat{N}\to H^2(G,\TT)$ is the only
obstruction for a central extension
$1\to N\to L\to G\to 1$
to be non-trivial. With a little bit of extra work one
can show that $H^2(G,\TT)$ is isomorphic to the nasty non-Hausdorff
group
 $\RR/(\ZZ+\theta\ZZ)$. In particular,
$G$ is not smooth, and the general structure theorem for
the group algebras of central extensions of $G$ as given
in Theorem \ref{thm-groupex} does not apply. However, some
weaker results can be deduced from \cite{en}.
\end{example}


\begin{thebibliography}{10}

\bibitem{andpas:hjm89}
Joel Anderson and William Paschke, \emph{The rotation algebra},
Houston J.\ Math. \textbf{15} (1989), 1--26.
\bibitem{blanchard2}
\'Etienne Blanchard, \emph{Tensor products of {$C(X)$}-algebras over {$C(X)$}},
   Ast{\'e}risque \textbf{232} (1995), 81--92.

\bibitem{en} Siegfried Echterhoff and Ryszard Nest,
\emph{The structure of the Brauer group and crossed products of
$C_0(X)$-linear group actions on $C_0(X,\K)$.}
Trans. Amer. Math. Soc. \textbf{353} (2001), 3685--3712.

\bibitem{ew3}
Siegfried Echterhoff and Dana~P. Williams, \emph{Crossed products by
   {$C_0(X)$}-actions}, J. Funct. Anal. \textbf{158} (1998), 113--151.

\bibitem{ew2}
\bysame, \emph{Locally inner actions on {$C_0(X)$}-algebras}, J. Operator
   Theory \textbf{45} (2001), 131--160.



\bibitem{HORR}
Steven Hurder, Dorte Olesen, Iain Raeburn, and Jonathan Rosenberg, \emph{The
   {C}onnes spectrum for actions of abelian groups on continuous-trace
   algebras}, Ergod. Th. \& Dynam. Sys. \textbf{6} (1986), 541--560.

\bibitem{moore1}
Calvin~C. Moore, \emph{Extensions and low dimensional cohomology theory of
   locally compact groups. {I}}, Trans. Amer. Math. Soc. \textbf{113}
   ({\noopsort{a}}1964), 40--63.

\bibitem{moore2}
\bysame, \emph{Extensions and low dimensional cohomology theory of locally
   compact groups. {II}}, Trans. Amer. Math. Soc. \textbf{113}
   ({\noopsort{b}}1964), 64--86.

\bibitem{moore3}
\bysame, \emph{Group extensions and cohomology for locally compact groups.
   {III}}, Trans. Amer. Math. Soc. \textbf{221} ({\noopsort{d}}1976), 1--33.

\bibitem{moore4}
\bysame, \emph{Group extensions and cohomology for locally compact groups.
   {IV}}, Trans. Amer. Math. Soc. \textbf{221} ({\noopsort{e}}1976), 34--58.

\bibitem{doir}
Dorte Olesen and Iain Raeburn, \emph{Pointwise unitary automorphism groups}, J.
   Funct. Anal. \textbf{93} (1990), 278--309.




\bibitem{judy-exp}
Judith~A. Packer, \emph{Transformation group {\cs}-algebras: a selective
   survey}, {\cs}-algebras: 1943--1993: A fifty year Celebration (Robert~S.
   Doran, ed.), Contemp. Math., vol. 169, Amer. Math. Soc., 1994, pp.~182--217.



\bibitem{judy96}
Judith~A. Packer, \emph{Moore cohomology and central twisted crossed product
   ${C}^ *$-algebras}, Canad. J. Math. \textbf{48} (1996), 159--174.
  
\bibitem{para1}
Judith~A. Packer and Iain Raeburn, \emph{Twisted crossed products of
   {\cs}-algebras}, Math. Proc. Camb. Phil. Soc. \textbf{106} (1989), 293--311.

\bibitem{para3}
Judith~A. Packer and Iain Raeburn, \emph{On the structure of
twisted group  {\cs}-algebras}, Trans. Amer. Math. Soc. \textbf{334} 
(1992), 685--718.


\bibitem{palais}
Richard~S. Palais, \emph{On the existence of slices for actions of non-compact
   {L}ie groups}, Ann. of Math. \textbf{73} (1961), 295--323.

\bibitem{quigg86}
John~C. Quigg, \emph{Duality for reduced twisted crossed products of
   $\cs$-algebras}, Indiana Univ. Math. J. \textbf{35} (1986), 549--571.
  
\bibitem{rr}
Iain Raeburn and Jonathan Rosenberg, \emph{Crossed products of continuous-trace
   {\cs}-algebras by smooth actions}, Trans. Amer. Math. Soc. \textbf{305}
   (1988), 1--45.

\bibitem{rw}
Iain Raeburn and Dana~P. Williams, \emph{Pull-backs of {\cs}-algebras and
   crossed products by certain diagonal actions}, Trans. Amer. Math. Soc.
   \textbf{287} (1985), 755--777.

\bibitem{90a}
\bysame, \emph{Moore cohomology, principal bundles, and actions of groups on
   {\cs}-algebras}, Indiana Univ. Math. J. \textbf{40} (1991), 707--740.

\bibitem{rw-book}
\bysame, \emph{Morita equivalence and continuous-trace {\cs}-algebras},
   Mathematical Surveys and Monographs, vol.~60, American Mathematical Society,
   Providence, RI, 1998.

\bibitem{ros2}
Jonathan Rosenberg, \emph{Some results on cohomology with {B}orel cochains,
   with applications to group actions on operator algebras}, Operator Theory:
   Advances and Applications \textbf{17} (1986), 301--330.

\bibitem{hasm1}
Harvey~A. Smith, \emph{Commutative twisted group algebras},
Trans. Amer. Math. Soc. 
   \textbf{197} (1974), 315--326.

\bibitem{hasm2}
\bysame, \emph{Characteristic principal bundles},
Trans. Amer. Math. Soc. 
   \textbf{211} (1975), 365--375.
   
 \bibitem{hasm3} \bysame, \emph{Central twisted group
     algebras}, Trans. Amer. Math. Soc.  \textbf{238} (1978),
   309--320.

\end{thebibliography}

\def\noopsort#1{}\def\cprime{$'$}
\providecommand{\bysame}{\leavevmode\hbox to3em{\hrulefill}\thinspace}
\providecommand{\MR}{\relax\ifhmode\unskip\space\fi MR }
\providecommand{\MRhref}[2]{%
   \href{http://www.ams.org/mathscinet-getitem?mr=#1}{#2}
}
\providecommand{\href}[2]{#2}

\end{document}